\documentclass[11pt]{amsart}

\usepackage{amsmath}
\usepackage{amsxtra}
\usepackage{amscd}
\usepackage{amsthm}
\usepackage{amsfonts}
\usepackage{amssymb}
\usepackage{eucal}

\newcommand{\slth}{\widehat{\mathfrak{sl}}_2}

\newcommand{\bra}[1]{\langle #1 |}        %bra
\newcommand{\ket}[1]{{| #1 \rangle}}      %ket
\newcommand{\Hom}{\mathop{\rm Hom}}
\newcommand{\cO}{\mathcal{O}}
\newcommand{\F}{\mathcal{F}}
\newcommand{\G}{\mathcal{G}}
\newcommand{\cH}{\mathcal{H}}
\newcommand{\cL}{\mathcal{L}}

\newcommand{\nn}{\nonumber}
\newcommand{\bea}{\begin{eqnarray}}
\newcommand{\ena}{\end{eqnarray}}
\newcommand{\be}{\begin{eqnarray*}}
\newcommand{\en}{\end{eqnarray*}}
\newcommand{\Vir}{{\it Vir}}
\newcommand{\cC}{\check{C}}
\newcommand{\yc}{\check{y}}

\newcommand{\C}{{\mathbb C}}
\newcommand{\Q}{{\mathbb Q}}
\newcommand{\Z}{{\mathbb Z}}

\numberwithin{equation}{section}
\newtheorem{thm}{Theorem}[section]
\newtheorem{prop}[thm]{Proposition}
\newtheorem{lem}[thm]{Lemma}

\newtheorem{rem}[thm]{Remark}

\begin{document} 

\title[Monomial basis for $M(p,p')$]
{A monomial basis for the Virasoro minimal series
$M(p,p')$~:~the case $1<p'/p<2$}

\author{B. Feigin, M. Jimbo, T. Miwa, E. Mukhin and Y. Takeyama}
\address{BF: Landau institute for Theoretical Physics, 
Chernogolovka,
142432, Russia}\email{feigin@feigin.mccme.ru}  
\address{MJ: Graduate School of Mathematical Sciences, 
The University of Tokyo, Tokyo 153-8914, Japan}
\email{jimbomic@ms.u-tokyo.ac.jp}
\address{TM: Department of Mathematics, 
Graduate School of Science,
Kyoto University, Kyoto 606-8502,
Japan}\email{tetsuji@math.kyoto-u.ac.jp}
\address{EM: Department of Mathematics,
Indiana University-Purdue University-Indianapolis,
402 N.Blackford St., LD 270,
Indianapolis, IN 46202}\email{mukhin@math.iupui.edu}
\address{YT: Institute of Mathematics, 
University of Tsukuba, Tsukuba, Ibaraki 305-8571,
Japan}\email{takeyama@math.tsukuba.ac.jp}

\date{\today}

\begin{abstract}
Quadratic relations of the intertwiners are given explicitly
in two cases of chiral conformal field theory, and monomial bases
of the representation spaces are constructed by using the Fourier components
of the intertwiners. The two cases are
the $(p,p')$-minimal series $M_{r,s}$ $(1\leq r\leq p-1,1\leq s\leq p'-1)$
for the Virasoro algebra where $1<p'/p<2$, and the level $k$ integrable
highest weight modules $V(\lambda)$ $(0\leq\lambda\leq k)$ for the affine Lie
algebra $\widehat{\mathfrak{sl}}_2$.
\end{abstract}
%%%%%%%%%%%%%%%%%%%%%%%%%%%%%%%%%%%%%
\maketitle

\setcounter{section}{0}
\setcounter{equation}{0}

%%%%%%%%%%%%%%%%%%%%%%%%%% ??? %%%%%%%%%%%%%%%%%%%%%%%%%%%%%%
\section{Introduction}
In this paper, we study quadratic relations satisfied
by intertwiners (or primary fields)
in chiral conformal field theory in two basic cases:
(I) the $(p,p')$-minimal series $M_{r,s}$ $(1\leq r\leq p-1,1\leq s\leq p'-1)$
for the Virasoro algebra where $1<p'/p<2$, and (II) the level $k$ integrable
highest weight modules $V(\lambda)$ $(0\leq\lambda\leq k)$ for the affine Lie
algebra $\widehat{\mathfrak{sl}}_2$. We consider special intertwiners:
the $(2,1)$-primary field acting as $M_{r,s}\rightarrow M_{r\pm1,s}$ for (I),
and the $\C^2$ intertwiner acting as
$V(\lambda)\rightarrow V(\lambda\pm1)$ for (II). We write explicitly the
quadratic relations for these intertwiners. Our aim is to construct
bases of representation spaces by using the intertwiners.
We construct vectors in representation spaces from the highest weight vectors
by the action of monomials in the Fourier components of intertwiners.
In Case (I), the space $\oplus_{r=1}^{p-1}M_{r,s}$ is
generated from the highest weight vector $|b(s),s\rangle\in M_{b(s),s}$ by
the action of the $(2,1)$-primary field. We choose $1\leq b(s)\leq p-1$
in such a way that the conformal dimension $\Delta_{r,s}$
of the space $M_{r,s}$ attains the minimum at $r=b(s)$. Similarly,
in Case (II), the space $\oplus_{\lambda=0}^kV(\lambda)$ is generated
from the highest weight vector $|0\rangle\in V(0)$ by the action of the
$\C^2$ intertwiner. We call these vectors monomials.
Using the quadratic relations, we rewrite an arbitrary
monomial as a linear combination of monomials satisfying certain
admissibility conditions. We then prove that admissible monomials are
linearly independent by computing the characters of representations.

There exist several constructions of monomial bases in terms of the chiral
conformal algebras (the Virasoro algebra, the affine Lie algebras, and so on),
instead of the intertwiners \cite{LW,LP,P,MP}.
For comparison we review one of such constructions for the Virasoro algebra.
In \cite{FF} a monomial base of the irreducible module
$M_{1,s}(\simeq M_{1,2k+3-s})$ $(1\leq s\leq k+1)$ in the $(2,2k+3)$-minimal
series is constructed. In addition to the Lie algebra relations,
which is quadratic, the Virasoro generators $L_n$ satisfy additional relations
of order $k+1$: for each $m\in\Z$ an element of the form
\begin{equation*}
\sum_{l=0}^k\sum_{n_1\leq\cdots\leq n_l\atop n_1+\cdots+n_l=m}
c^{(m)}_{n_1,\ldots,n_l}L_{n_1}\cdots L_{n_l}
\end{equation*}
acts as zero on all modules $M_{1,s}$. These relations are obtained from a
single relation for the Virasoro current. Using the relations one can reduce
any monomial of the form
\begin{equation*}
L_{-n_1}\cdots L_{-n_N}|1,s\rangle\quad(n_1\geq\cdots\geq n_N>0)
\end{equation*}
to those which satisfy the difference two conditions
\begin{equation}\label{DIFTWO}
n_i-n_{i+k}\geq 2\quad(1\leq i\leq N-k).
\end{equation}

In the case of intertwiners, the operators act from one irreducible module
to another. Therefore, in order to label monomials we need to specify
a sequence of modules which appear in the composition of intertwiners.
We call it a path. Since different modules have different conformal dimensions,
a natural parametrization of the Fourier components is by rational numbers
instead of integers.  In Case (I), monomials are of the form
\begin{equation*}
\phi^{(r_0,r_1)}_{-n_1}\cdots \phi^{(r_{N-1},r_N)}_{-n_N}|b(s),s\rangle
\end{equation*}
where $r_N=b(s),r_i=r_{i+1}\pm1$ and
$n_i\in\Delta_{r_{i-1},s}-\Delta_{r_i,s}+\Z$.
The operator $\phi^{(r,r')}_n$ maps $M_{r',s}\rightarrow M_{r,s}$.

The difference two conditions change to
\begin{equation*}
n_i-n_{i+1}\geq w(r_i,r_{i+1},r_{i+2})\quad(1\leq i\leq N-1)
\end{equation*}
where $w(r,r',r'')$ $(r=r'\pm1,r'=r''\pm1)$ are rational numbers.
Namely, the distance $k$ and the gap $2$ in (\ref{DIFTWO}) change to
the distance $1$ and the gap $w$. The gap is dependent on the path.

Let us illustrate the gap condition in the simplest case $w(r+2,r+1,r)=t/2$.
The quadratic relation reads
\begin{eqnarray*}
&&z_1^{-t/2}(1-z_2/z_1)^{-t/2}\phi^{(r+2,r+1)}(z_1)\phi^{(r+1,r)}(z_2)\\
&&=z_2^{-t/2}(1-z_1/z_2)^{-t/2}\phi^{(r+2,r+1)}(z_2)\phi^{(r+1,r)}(z_1),
\end{eqnarray*}
where we have the expansion
\begin{equation*}
\phi^{(r,r')}(z)=\sum_{n\in\Z+\Delta_{r',s}-\Delta_{r,s}}
\phi_nz^{-n-\Delta_{2,1}}.
\end{equation*}
Let $(1-z)^{-t/2}=\sum_{n\geq0}c_nz^n$ be the Taylor expansion.
In terms of the Fourier coefficients we have the following relations
\begin{eqnarray*}
&&\phi^{(r+2,r+1)}_m\phi^{(r+1,r)}_n
+c_1\phi^{(r+2,r+1)}_{m-1}\phi^{(r+1,r)}_{n+1}+\cdots\\
&&=\phi^{(r+2,r+1)}_{n-t/2}\phi^{(r+1,r)}_{m+t/2}
+c_1\phi^{(r+2,r+1)}_{n-t/2-1}\phi^{(r+1,r)}_{m+t/2+1}+\cdots.
\end{eqnarray*}
If $n-m<t/2$, we can reduce a monomial $\phi^{(r+2,r+1)}_m\phi^{(r+1,r)}_n$
by this relation to a linear combination of monomials
$\phi^{(r+2,r+1)}_{m'}\phi^{(r+1,r)}_{n'}$ where $n'-m'>n-m$.
The case $w(r,r\pm1,r)$ is more involved. Instead of the power function
$(1-z_2/z_1)^{-t/2}$ we need a hypergeometric series. The details
are given in the main text.

The combinatorial structure of monomial bases for Case (II) is similar.
In this case, the intertwiner has two components corresponding to
two weight vectors in the $\mathfrak{sl}_2$ module $\C^2$.
Therefore, we need another kind of paths which describe
the composition sequence of these components.

In \cite{FM} and \cite{FJM}, respectively, a monomial base in terms of the
$\C^{k+1}$ intertwiner for the level $k$ integrable highest weight
$\widehat{\mathfrak{sl}}_2$ modules, and one in terms of the $(2,1)$-primary
field in the $(3,p')$ Virasoro minimal series, respectively, is constructed, 
To be precise, monomial bases constructed in these papers are not bases
for the modules but ones for direct sums of the tensor products of the modules
with certain bosonic Fock spaces. The intertwiners are modified by
bosonic vertex operators so that they constitute a new VOA.

Finally, we mention \cite{H}, in which, in the framework of VOA,
an existence of an infinite set of quadratic relations is shown.
In the present paper, in Cases (I) and (II), we give explicit forms
of quadratic relations. We use finitely many components of intertwiners,
only one in Case (I) and two in Case (II), and write finitely many quadratic
relations. Each relation consists of infinitely many relations for the
Fourier components of the intertwiners. We also prove that the left ideal
which annihilates the highest weight vector is generated by these relations
(and the highest weight conditions) by showing that monomial bases are
obtained by reduction using these quadratic relations.

The plan of the paper is as follows. In Section 2 we derive the quadratic
relations of the intertwiners for Case (I).
In Section 3 we construct monomial bases for Case (I). Construction of the
quadratic relations and the monomial bases for Case (II)
is given in Section 4.

%%%%%%%%%%%%%%%%%%%%%%%%%% ??? %%%%%%%%%%%%%%%%%%%%%%%%%%%%%%

\section{Quadratic relations 
for the $(2,1)$-primary field}\label{sec:exchange}

\subsection{The $(2,1)$-primary field in the minimal series}
\label{subsec:notation}
In this subsection we introduce our notation, 
and summarize some basic facts concerning 
representations of the Virasoro algebra  
which will be used in subsequent sections. 
For these and related formulas, we find 
the textbook \cite{DiFMS} to be a useful reference. 

Let $\Vir$ be the Virasoro algebra with the standard $\C$-basis 
$\{L_n\}_{n\in\Z}$ and $c$, satisfying
\be
[L_m,L_n]=(m-n)L_{m+n}+\frac{c}{12}m(m^2-1)\delta_{m+n,0},
\quad [c,L_n]=0.
\en
Fix a pair $(p,p')$ of relatively prime positive integers.  
We set 
\be
t=\frac{p'}{p}. 
\en 
Though we will consider only the case $1<t<2$ later on,  
we do not make this assumption in Section \ref{sec:exchange}.
For each $(p,p')$, there exists a family of 
irreducible $\Vir$-modules 
$M_{r,s}=M_{r,s}(p,p')$ ($1\le r\le p-1,~1\le s\le p'-1$) 
on which $c$ acts as the scalar
\be
c_{p,p'}%=1-\frac{6(p-p')^2}{4pp'}
=13-6\left(t+\frac{1}{t}\right). 
\en
We use the notation $\ket{v}$ to represent an element of $M_{r,s}$. 
The module $M_{r,s}$ is $\Q$-graded with respect to the degree operator 
$L_0$, 
\be
&&M_{r,s}=\oplus_{d\in \Z_{\ge 0}+\Delta_{r,s}}(M_{r,s})_d,
\\
&&
(M_{r,s})_d=
\{\ket{v}\in M_{r,s}\mid L_0\ket{v}=d\ket{v}\},
\en
where $\dim_\C (M_{r,s})_{\Delta_{r,s}}=1$ and
\be
&&\Delta_{r,s}=\frac{(rt-s)^2-(t-1)^2}{4t}. 
\en
We fix a generator of $(M_{r,s})_{\Delta_{r,s}}$ and denote it by $\ket
{r,s}$. 
We have $L_n(M_{r,s})_d\subset (M_{r,s})_{d-n}$, so that 
$L_n\ket{r,s}=0$ ($n>0$), 
and $L_{-n}$ ($n>0$) are the creation operators. 

Consider the right $\Vir$-module 
$M_{r,s}^*=\oplus_{d\in \Z_{\ge 0}+\Delta_{r,s}}(M_{r,s})^*_d$, 
$(M_{r,s})^*_d=\Hom_\C\bigl((M_{r,s})_d,\C\bigr)$. 
For $\bra{u}\in (M_{r,s})_d^*$ and $\ket{v}\in(M_{r,s})_d$, 
we set $\deg\bra{u}=d$, $\deg\ket{v}=d$, and  
write the dual coupling $M_{r,s}^*\times M_{r,s}\rightarrow \C$ 
as $(\bra{u},\ket{v})\mapsto 
\langle u|v \rangle$. 
We denote by $\bra{r,s}\in (M_{r,s})_{\Delta_{r,s}}^*$ 
the element such that $\langle{r,s}|{r,s}\rangle=1$. 

In conformal field theory, 
the notion of primary fields plays a key role.  
The $(k,l)$-primary field is a collection of generating series 
\be
&&\phi^{(r',s';r,s)}_{k,l}(z)=
\sum_{n\in\Z+\Delta_{r,s}-\Delta_{r',s'}}
(\phi^{(r',s';r,s)}_{k,l})_nz^{-n-\Delta_{k,l}},
\en
whose coefficients are linear operators 
$(\phi^{(r',s';r,s)}_{k,l})_n:M_{r,s}\rightarrow M_{r',s'}$. 
Up to a scalar multiple, 
it is characterized by the following properties.
\bea
&&[L_n,\phi^{(r',s';r,s)}_{k,l}(z)]=
z^n\bigl(z\partial+(n+1)\Delta_{k,l}\bigr)
\phi^{(r',s';r,s)}_{k,l}(z),
\label{intertwiner}
\\
&&(\phi^{(r',s';r,s)}_{k,l})_n(M_{r,s})_d
\subset (M_{r',s'})_{d-n}.\notag
\label{homogeneity}
%\\
%&&(\phi^{(r',s';r,s)}_{k,l})_{\Delta_{r,s}-\Delta_{r',s'}}
%\ket{r,s}=\ket{r',s'}.
%\label{normalization_of_primary}
\ena
Here and after, we set $\partial=d/dz$. 

In this paper we will consider only the $(2,1)$-primary field. 
It exists for $(r',s')=(r\pm 1,s)$ with 
$1\le r,r\pm1\le p-1$, $1\le s\le p'-1$. 
Normally we suppress the index $(r,s)$ and write it as 
\bea
\phi^\pm(z)=\phi_{2,1}^{(r\pm 1,s;r,s)}(z)
=\sum_{n\in\Z+\Delta_{r,s}-\Delta_{r\pm 1,s}}
\phi_n^\pm\, z^{-n-\Delta_{2,1}},
\label{phi21}
\ena
where $\Delta_{2,1}=(3t-2)/4$.
We choose the normalization
\footnote{Our normalization \eqref{normalization_of_primary} 
is different from the one commonly used 
in conformal field theory, 
where one demands 
$\bra{r,s}\phi^{\sigma}(z_1)\phi^{-\sigma}(z_2)\ket{r,s}
=(z_1-z_2)^{-2\Delta_{2,1}}+\cdots$ as $z_1\to z_2$.}
\bea
\phi_{\Delta_{r,s}-\Delta_{r\pm 1,s}}^\pm\ket{r,s}
=\ket{r\pm 1,s}.
\label{normalization_of_primary}
\ena
It is known that the highest-to-highest matrix element 
\be
H(z_1,z_2)=\bra{r',s}\phi^{\sigma_1}(z_1)\phi^{\sigma_2}(z_2)\ket{r,s}
\quad (r'=r+\sigma_1+\sigma_2)
\en 
satisfies the second order linear differential equation
\be
\left(\frac{1}{t}\partial^2
+\Bigl(\frac{1}{z}+\frac{1}{z-1}\Bigr)\partial
{}-\frac{\Delta_{r,s}}{z^2}-\frac{\Delta_{2,1}}{(z-1)^2}
+\frac{2\Delta_{2,1}+\Delta_{r,s}-\Delta_{r',s}}{z(z-1)}\right)
H(1,z)=0.
\en
The existence, and the above 
differential equation, 
are the only information we will need about $\phi^\pm(z)$.

Together with the homogeneity
$H(kz_1,kz_2)=
k^{\Delta_{r',s}-\Delta_{r,s}-2\Delta_{2,1}}H(z_1,z_2)$ 
and the normalization \eqref{normalization_of_primary}, 
$H(z_1,z_2)$ are determined as follows. 
\be
&&
\bra{r\pm 2,s}\phi^\pm(z_1)\phi^\pm(z_2)\ket{r,s}
=(z_1z_2)^{-(t-1)/2\pm(rt-s)/2}(z_1-z_2)^{t/2},\\
&&
\bra{r,s}\phi^\mp(z_1)\phi^\pm(z_2)\ket{r,s}
=(z_1z_2)^{-\Delta_{2,1}}y_{r,s}^{\pm}(z_2/z_1),
\en
where
\bea\label{Y}
&&
y_{r,s}^{\pm}(z)=
z^{t/4\pm(rt-s)/2}(1-z)^{1-3t/2}
F\bigl(1-t,1-t\pm(rt-s),1\pm(rt-s);z\bigr)
\ena
and $F(a,b,c;z)$ denotes the hypergeometric function.
In the particular case $r=1$ or $r=p-1$, 
there are only the following possibilities:
\be
\mbox{$r=1$}&:& 
y_{1,s}^+(z)=z^{3t/4-s/2}(1-z)^{1-3t/2}F(1-t,1-s,1-s+t;z),
\\
\mbox{$r=p-1$}&:&
y_{p-1,s}^-(z)=
z^{3t/4-(p'-s)/2}(1-z)^{1-3t/2}F(1-t,1-p'+s,1-p'+s+t;z).
\en
The relevant hypergeometric functions are reciprocal polynomials, i.e., 
\bea
F(a,-n,1-a-n;z^{-1})=z^{-n}F(a,-n,1-a-n;z)\qquad (n\in\Z_{\ge 0}).
\label{reciprocal}
\ena

{}From these results 
and the intertwining property \eqref{intertwiner}, 
it follows that 
the general matrix elements $\bra{u}\phi^{\sigma_1}(z_1)\phi^{\sigma_2}
(z_2)\ket{v}$ 
($\bra{u}\in M_{r',s}^*,\ket{v}\in M_{r,s}$) 
are Laurent series convergent in the domain $|z_1|>|z_2|$, 
multiplied by an overall rational power of $z_1,z_2$. 

\subsection{A bilinear relation for hypergeometric functions}
The hypergeometric equation for $y_{r,s}^{\pm}(z)$ is invariant
under the substitution $z\to z^{-1}$. In this subsection, we derive
an identity for solutions to such an equation.

Let us consider a Fuchsian linear differential equation with 
regular singularities at $0,1,\infty$, 
which is also invariant under the substitution $z\to z^{-1}$. 
It takes the general form 
\be
\frac{d^2y}{dz^2}+
\left(\frac{1-\lambda_+-\lambda_-}{z}+\frac{1-\mu-\nu}{z-1}
\right)\frac{dy}{dz}
+\left(\frac{\lambda_+\lambda_-}{z^2}
+\frac{\mu\nu}{z(z-1)^2}\right)y=0,
\label{HGE}
\en
where $\lambda_\pm,\mu,\nu\in\C$ are 
parameters satisfying the relation 
\be
2(\lambda_+ + \lambda_-)+\mu+\nu=1.
\en
The corresponding Riemann scheme is
\bea
\left\{
\begin{matrix}
0 & 1& \infty \\
\lambda_+&\mu&\lambda_+\\
\lambda_-&\nu&\lambda_-\\
\end{matrix}
\right\}. 
\label{Riemann}
\ena
We assume that $\lambda_+-\lambda_-,\mu-\nu\not\in\Z$. 
Then a basis of solutions is given by 
\be
&&y^\pm(z)=C^\pm (-z)^{\lambda_\pm}(1-z)^\mu
F\bigl(\lambda_++\lambda_-+\mu,
2\lambda_\pm+\mu,
1\pm(\lambda_+-\lambda_-);z\bigr)
\en
with $C^\pm\neq 0$. 
Fixing the branch by ${\rm arg}\,(-z)=0$ for $z<0$, 
we have the transformation law 
\be
y^{\tau}(z^{-1})=\sum_{\sigma=\pm} 
y^\sigma(z)B_\sigma^{\tau},
\en
where the connection matrix 
$B=(B_\sigma^\tau)_{\sigma,\tau=\pm}$ is given by 
\bea
B_\pm^\pm=\frac{\sin\pi(\lambda_++\lambda_-+\mu)}
{\sin\pi(\lambda_\mp-\lambda_\pm)},
\quad
B_\mp^\pm=\frac{
\Gamma(\lambda_\pm-\lambda_\mp)\Gamma(1+\lambda_\pm-\lambda_\mp)}
{\Gamma(2\lambda_\pm+\mu)\Gamma(1-2\lambda_\mp-\mu)}
\frac{C^\pm}{C^\mp}.
\label{braid}
\ena

Consider now the Riemann scheme 
\bea
\left\{
\begin{matrix}
0 & 1& \infty \\
{}-\lambda_+&-\mu&-\lambda_+\\
{}-\lambda_-&2-\nu&-\lambda_-\\
\end{matrix}
\right\}.
\label{dualRiemann}
\ena
Let 
\bea
&&\yc^\pm(z)=\cC^\pm (-z)^{-\lambda_\pm}(1-z)^{-\mu}
F\bigl(-\lambda_+-\lambda_- -\mu,
{}-2\lambda_\pm-\mu,
1\mp(\lambda_+-\lambda_-);z\bigr)
\label{ycheck}
\ena
be a basis of solutions of 
the corresponding differential equation. 
{}From \eqref{braid} it follows that, if
\bea
{C^+\cC^+}
=-\frac{2\lambda_-+\mu}{2\lambda_++\mu}\,{C^-\cC^-},
\label{cpcp}
\ena
then the connection matrix
associated with \eqref{ycheck} is given by ${}^tB^{-1}$. 
We take
\bea
\cC^\pm C^\pm=\frac{2\lambda_\mp+\mu}{2(\lambda_\mp-\lambda_\pm)},
\label{CCc}
\ena
so that \eqref{cpcp} and $C^+\cC^++C^-\cC^-=1$ hold. 

\begin{lem}\label{lem:bilinear}
With the choice \eqref{CCc}, the following identities hold.
\bea
&&y^+(z)\yc^+(z)+y^-(z)\yc^-(z)=1,
\label{bilin1}\\
&&y^+(z)(z\partial\yc^+)(z)+y^-(z)(z\partial\yc^-)(z)=\frac{\mu}{2}\,
\frac{1+z}{1-z}.
\label{bilin2}
\ena
\end{lem}
\begin{proof}
Denote by $\varphi(z)$, $\tilde{\varphi}(z)$ 
the left hand sides of \eqref{bilin1}, \eqref{bilin2}, 
respectively. 
They are single valued and holomorphic at $z=0$. 
{}From the relations of the connection matrices mentioned above, 
we have $\varphi(z^{-1})=\varphi(z)$,  
$\tilde{\varphi}(z^{-1})=-\tilde{\varphi}(z)$.  
Therefore they are single valued also at $z=\infty$, and 
hence on $\mathbb{P}^1$. 
{}From the Riemann schemes \eqref{Riemann}, \eqref{dualRiemann}, 
we see that $\varphi(z)$ is regular at $z=1$, 
while $\tilde{\varphi}(z)$ has at most a simple pole there.
The lemma follows from these facts.
\end{proof}

\subsection{Exchange relations}
The aim of this subsection is to derive the 
quadratic exchange relations for the $(2,1)$-primary field. 

First let us introduce several functions 
$f^\sigma_{r,s},g^\sigma_{r,s},h^\sigma_{r,s}$ 
which enter these relations. 
They are of the form $z_1^\alpha z_2^\beta \psi(z_2/z_1)$, 
where $\alpha,\beta\in\Q$ and $\psi(z)$ is a power series convergent in 
$|z|<1$. 

Suppose $2\le r\le p-2$. 
Set
\be
&&\lambda_\pm=\Delta_{r\pm 1,s}-\Delta_{r,s}=\frac{t}{4}\pm\frac{rt-s}
{2},
\\
&&\mu=-2\Delta_{2,1}=1-\frac{3t}{2},\quad
\quad
\nu=\Delta_{3,1}-2\Delta_{2,1}=\frac{t}{2}, 
\en
and $C^{\pm}=e^{\pi i\lambda_\pm}$. 

With this choice of parameters, $y^\pm(z)$ becomes $y^\pm_{r,s}(z)$ 
in \eqref{Y}. 
The corresponding functions $\yc^\pm(z)$ defined by \eqref{ycheck} and 
\eqref{CCc} will be denoted by $\yc^\pm_{r,s}(z)$.

Define
\bea
&&f_{r,s}^\pm(z_1,z_2)=(z_1z_2)^{\Delta_{2,1}}
z_1^{-1}(1-z_2/z_1)^{-1}\yc_{r,s}^\pm(z_2/z_1),
\label{fsigma} \notag
\\
&&g_{r,s}^\pm(z_1,z_2)=(z_1z_2)^{\Delta_{2,1}}
(z\partial\yc^{\pm}_{r,s})(z_2/z_1).
\label{gsigma} \notag
\ena
In the case $r=1$ or $p-1$, we define 
$f_{r,s}^\sigma(z_1,z_2)=g_{r,s}^\sigma(z_1,z_2)=0$ 
except for the following ones. 
\bea
&&
f_{1,s}^+(z_1,z_2)=
K_{1-t,s-1}
z_1^{3(t-1)/2+s/2-1}z_2^{(s-1)/2}(1-z_2/z_1)^{-2+3t/2},
\label{1.18}
\\
&&f_{p-1,s}^-(z_1,z_2)=
K_{1-t,p'-s-1}
z_1^{3(t-1)/2+(p'-s)/2-1}z_2^{(p'-s-1)/2}(1-z_2/z_1)^{-2+3t/2}.
\label{1.19} \notag
\ena
Here $K_{a,n}=F(a,-n,1-a-n;1)^{-1}=\prod_{j=0}^{n-1}(a+j)/(2a+j)$.
Define also
\bea
&&
h_{r,s}^\pm(z_1,z_2)=(z_1z_2)^{(t-1)/2\mp(rt-s)/2}
z_1^{-t/2}(1-z_2/z_1)^{-t/2}
\label{1.20}
\ena
for all $r,s$. 

Consider the formal series
\be
&&
\F_{r,s}(z_1,z_2):=
\sum_{\sigma=\pm}
f_{r,s}^\sigma(z_1,z_2)\phi^{-\sigma}(z_1)\phi^{\sigma}(z_2),
\\
&&\G_{r,s}(z_1,z_2):=
\sum_{\sigma=\pm}
g_{r,s}^\sigma(z_1,z_2)\phi^{-\sigma}(z_1)\phi^{\sigma}(z_2),
\\
&&\cH^\sigma_{r,s}(z_1,z_2):=
h_{r,s}^\sigma(z_1,z_2)\phi^\sigma(z_1)\phi^\sigma(z_2).
\en
These series have the form $\sum_{m,n\in\Z}\cO_{m,n}z_1^mz_2^n$, 
where each coefficient $\cO_{m,n}$ is a well defined linear operator. 

\begin{lem}\label{lem:exchange-vac}
The following identities hold:
\bea
&&\bra{r,s}\F_{r,s}(z_1,z_2)\ket{r,s}
=\frac{1}{z_1-z_2}\,P_{r,s}(z_1,z_2),
\label{ex-vacf}
\\
&&\bra{r,s}\G_{r,s}(z_1,z_2)\ket{r,s}
=
(-\Delta_{2,1})\frac{z_1+z_2}{z_1-z_2}\qquad (2\le r\le p-2),
\label{ex-vacg}
\\
&&
\bra{r+2\sigma,s}\cH^\sigma(z_1,z_2)\ket{r,s}= 1,
\label{ex-vach}
\ena
where $P_{r,s}(z_1,z_2)$ is a homogeneous symmetric polynomial 
satisfying
\be
P_{r,s}(z,z)=
\begin{cases}
1& (2\le r\le p-2),\\
z^{s-1}& (r=1),\\
z^{p'-s-1}& (r=p-1).\\
\end{cases}
\en
\end{lem}
\begin{proof}
Using Lemma \ref{lem:bilinear}, formulas 
\eqref{ex-vacf}, \eqref{ex-vacg} for $2\le r\le p-2$
and \eqref{ex-vach} are easily verified with $P_{r,s}(z_1,z_2)=1$. 
Formula \eqref{ex-vacf} for $r=1$ follows from 
\be
P_{1,s}(z_1,z_2)
=z_1^{s-1}
\frac{F(1-t,1-s,1-s+t;z_2/z_1)}{F(1-t,1-s,1-s
+t;1)}
\en
and \eqref{reciprocal}. The case $r=p-1$ is similar. 
\end{proof}

Denote by $\C[z_1^{\pm 1},z_2^{\pm 1}]^{S_2}$ 
the space of symmetric Laurent polynomials in $z_1,z_2$. 
\begin{lem}\label{lem:FGH}
Set $\cL_n=z_1^{n+1}\partial_1
+z_2^{n+1}\partial_2+(n+1)(z_1^n+z_2^n)\Delta_{2,1}$,
$\partial_j=\partial/\partial z_j$. 
Then there exist 
$a_n,b_n,c_n,d_n,e_n\in \C[z_1^{\pm 1},z_2^{\pm 1}]^{S_2}$ 
such that 
\bea
&&
[L_n,\F_{r,s}(z_1,z_2)]=\bigl(\cL_n+a_n(z_1,z_2)\bigr)
\F_{r,s}(z_1,z_2)+b_n(z_1,z_2)\G_{r,s}(z_1,z_2),
\label{LF}
\\
&&
[L_n,\G_{r,s}(z_1,z_2)]=c_n(z_1,z_2)\F_{r,s}(z_1,z_2)+
\bigl(\cL_n+d_n(z_1,z_2)\bigr)\G_{r,s}(z_1,z_2),
\label{LG}
\\
&&
[L_n,\cH^\sigma_{r,s}(z_1,z_2)]=
\bigl(\cL_n+e_n(z_1,z_2)\bigr)\cH^\sigma_{r,s}(z_1,z_2).
\label{LH}
\ena
\end{lem}
\begin{proof}
Let us show \ref{LF} assuming $1<r<p-1$. Using 
\[
[L_n,\phi^{-\sigma}(z_1)\phi^{\sigma}(z_2)]
=\cL_n\phi^{-\sigma}(z_1)\phi^{\sigma}(z_2),
\]
we obtain
\be
&&
[L_n,\F_{r,s}(z_1,z_2)]-\cL_n\F_{r,s}(z_1,z_2)
\\
&&
=
{}-\sum_{\sigma=\pm}
(z_1^{n+1}\partial_1+z_2^{n+1}\partial_2)f_{r,s}^\sigma(z_1,z_2)
\cdot \phi^{-\sigma}(z_1)\phi^{\sigma}(z_2)
\\
&&=a_n(z_1,z_2)\F_{r,s}(z_1,z_2)+b_n(z_1,z_2)\G_{r,s}(z_1,z_2),
\en
where
\be
&&a_n(z_1,z_2)=
\frac{z_1^{n+1}-z_2^{n+1}}{z_1-z_2}-\Delta_{2,1}(z_1^n+z_2^n),
\\
&&b_n(z_1,z_2)=
\frac{z_1^n-z_2^n}{z_1-z_2}.
\en
Hence we have \eqref{LF}.
For $r=1$ or $p-1$, we find 
\be
[L_n,\F_{r,s}(z_1,z_2)]-\cL_n\F_{r,s}(z_1,z_2)
=a'_n(z_1,z_2)\F_{r,s}(z_1,z_2),
\en
where
\be
a'_n(z_1,z_2)=-\frac{s'-1}{2}(z_1^n+z_2^n)-\left(2\Delta_{2,1}-1\right)
\frac{z_1^{n+1}-z_2^{n+1}}{z_1-z_2}
\en 
and $s'=s$ ($r=1$) or $p'-s$ ($r=p-1$). 

Likewise we have
\be
&&[L_n,\G_{r,s}(z_1,z_2)]-\cL_n\G_{r,s}(z_1,z_2)
\\
&&=
{}-\Delta_{2,1}(z_1^n+z_2^n)\G_{r,s}(z_1,z_2)
\\
&&
+(z_1z_2)^{\Delta_{2,1}}(z_1^n-z_2^n)
\sum_{\sigma=\pm}
\left((z\partial)^2\yc_{r,s}^\sigma\right)(z_2/z_1)
\cdot \phi^{-\sigma}(z_1)\phi^{\sigma}(z_2).
\en
With the substitution
\be
(z\partial)^2\yc_{r,s}
=\frac{t}{2} \frac{z+1}{z-1}
z\partial \yc_{r,s}
{}-\left(\frac{t^2}{16}-\frac{(rt-s)^2}{4}-
(\frac{3t}{2}-1)(\frac{t}{2}-2)\frac{z}{(z-1)^2}
\right)\yc_{r,s}, 
\en
the right hand side becomes 
\be
c_n(z_1,z_2)\F_{r,s}(z_1,z_2)+d_n(z_1,z_2)\G_{r,s}(z_1,z_2),
\en
with
\be
&&c_n(z_1,z_2)=-\frac{z_1^n-z_2^n}{z_1-z_2}
\left(\Bigl(\frac{t^2}{16}-\frac{(rt-s)^2}{4}\Bigr)(z_1-z_2)^2
{}-(\frac{3t}{2}-1)(\frac{t}{2}-2)z_1z_2\right),
\\
&&d_n(z_1,z_2)=
-\left(\frac{t}{2}\frac{z_1^n-z_2^n}{z_1-z_2}(z_1+z_2)
+\Delta_{2,1}(z_1^n+z_2^n)\right).
\en
This proves \eqref{LG}.

The case \eqref{LH} is similar. 
\end{proof}

\begin{prop}\label{QUADREL}
Set $\delta(z_2/z_1)=\sum_{n\in\Z}z_1^{-n}z_2^n$.
We have the identities of formal series
\bea
&&
\sum_{\sigma=\pm}
f_{r,s}^\sigma(z_1,z_2)\phi^{-\sigma}(z_1)\phi^{\sigma}(z_2)
+\sum_{\sigma=\pm}
f_{r,s}^\sigma(z_2,z_1)\phi^{-\sigma}(z_2)\phi^{\sigma}(z_1)
\label{1.27}\\
&&\qquad\qquad
=z_1^{-1}\delta(z_2/z_1)\times
\begin{cases}
1& (2\le r\le p-2),\\
z_1^{s-1}& (r=1),\\
z_1^{p'-s-1}& (r=p-1),\\
\end{cases}
\nn
\\
&&
\sum_{\sigma=\pm}
g_{r,s}^\sigma(z_1,z_2)\phi^{-\sigma}(z_1)\phi^{\sigma}(z_2)
+\sum_{\sigma=\pm}
g_{r,s}^\sigma(z_2,z_1)\phi^{-\sigma}(z_2)\phi^{\sigma}(z_1)
\label{1.28}
\\
&&\qquad\qquad
=(-2\Delta_{2,1})\delta(z_2/z_1)
\qquad (2\le r\le p-2), \notag
\\
&&
h_{r,s}^\sigma(z_1,z_2)\phi^\sigma(z_1)\phi^\sigma(z_2)
=h_{r,s}^\sigma(z_2,z_1)\phi^\sigma(z_2)\phi^\sigma(z_1).
\label{1.29}
\ena
\end{prop}
\begin{proof}
We show that   
for any $\bra{u}\in M^*_{r',s}$, $\ket{v}\in M_{r,s}$ we have
modulo $(z_1-z_2)\C[z_1^{\pm 1},z_2^{\pm 1}]^{S_2}$ 
\bea
&&\bra{u}\F_{r,s}(z_1,z_2)\ket{v}
\equiv \frac{1}{z_1-z_2}\langle u|v\rangle
\,P_{r,s}(z_1,z_2),
\label{FF}
\\
&&\bra{u}\G_{r,s}(z_1,z_2)\ket{v}
\equiv (-\Delta_{2,1})
\frac{z_1+z_2}{z_1-z_2}\langle u|v\rangle
\qquad (2\le r\le p-2),
\label{GG}
\ena
where $P_{r,s}(z_1,z_2)$ is as in Lemma \ref{lem:exchange-vac}. 

We also show that modulo $\C[z_1^{\pm 1},z_2^{\pm 1}]^{S_2}$
\bea
\bra{u}\cH^\sigma_{r,s}(z_1,z_2)\ket{v}\equiv 0.
\label{HH}
\ena

The left hand sides are convergent in the domain $|z_1|>|z_2|$. 
Noting that 
\be
\frac{1}{z_1-z_2}\,P_{r,s}(z_1,z_2)+(z_1\leftrightarrow z_2)
=z_1^{-1}\delta(z_2/z_1)\times P_{r,s}(z_1,z_1),
\en
we obtain the proposition.

By Lemma \ref{lem:exchange-vac}, 
\eqref{FF}--\eqref{HH}
hold true 
in the case $\bra{u}=\bra{r',s}, \ket{v}=\ket{r,s}$.
By induction, suppose they are true for $\bra{u'},\ket{v'}$ with 
$\deg\bra{u'}+\deg\ket{v'}<\deg\bra{u}+\deg\ket{v}$.   
We may assume either $\bra{u}=\bra{u'}L_n$ or $\ket{v}=L_{-n}\ket{v'}$ 
with some $\bra{u'}$, $\ket{v'}$ and $n>0$. 
{}From Lemma \ref{lem:FGH}, we have in the first case
with $1<r<p-1$, 
\be
&&\bra{u'}L_{n}\F_{r,s}(z_1,z_2)\ket{v}
\\
&&\equiv \bra{u'}\F_{r,s}(z_1,z_2)L_{n}\ket{v} + 
(\cL_n+a_n(z_1,z_2))\bra{u'}\F_{r,s}(z_1,z_2)\ket{v}
+b_n(z_1,z_2)\bra{u'}\G_{r,s}(z_1,z_2)\ket{v}
\\
&&\equiv \bra{u'}L_{n}\ket{v}\frac{1}{z_1-z_2}\,P_{r,s}(z_1,z_2)
\\
&&
+\left((\cL_n+a_n(z_1,z_2))\frac{1}{z_1-z_2}
+b_n(z_1,z_2)(-\Delta_{2,1})\frac{z_1+z_2}{z_1-z_2}
\right)\langle{u'}|{v}\rangle.
\en
In the last line we used the induction hypothesis.
The second term in the right hand side is a Laurent polynomial.
Indeed, using the expressions for $a_n,b_n$ given 
in the proof of Lemma \ref{lem:FGH}, we find 
\be
&&(z_1-z_2)\left(
(\cL_n+a_n(z_1,z_2))\frac{1}{z_1-z_2}
+b_n(z_1,z_2)(-\Delta_{2,1})\frac{z_1+z_2}{z_1-z_2}\right)
\\
&&=\Delta_{2,1}
\left(n(z_1^n+z_2^n)-\frac{z_1^n-z_2^n}{z_1-z_2}(z_1+z_2)
\right).
\en
Since it vanishes at $z_2=z_1$, the assertion follows. 

In the same way, in the case $r=1$ or $p-1$, 
we are to check that 
$(\cL_n+a'_n(z_1,z_2))\frac{P_{r,s}(z_1,z_2)}{z_1-z_2}$
is a Laurent polynomial, where $a'_n$ is as in 
the proof of Lemma \ref{lem:FGH}. 
This can be verified by noting that 
$(z_1^{n+1}\partial_1+z_2^{n+1}\partial_2)P_{r,s}(z_1,z_1)
=(s'-1)z_1^{n+s'-1}$. 

The other cases can be proved in a similar manner.
\end{proof}

%\section{Monomial basis}
\section{Monomial bases in terms of the primary field $\phi_{2,1}$}
We fix coprime integers $p,p'$ ($t=p'/p$) satisfying
\begin{equation}
1<t<2,\label{<2}
\end{equation}
and consider the representations $M_{r,s}$ ($1\leq r\leq p-1,1\leq s\leq p'-1$)
of the Virasoro algebra in the $(p,p')$ minimal series. We construct a
monomial basis of $M_{r,s}$ by using the $(2,1)$-primary field.

\subsection{Spanning set of monomials}
In this subsection we construct a spanning set of vectors for each $M_{r,s}$.
In the next subsection, we prove that it constitutes a basis.

For each $1\leq s\leq p'-1$ we define $b(s)$ by the condition
that the conformal dimension $\Delta_{r,s}$ takes the minimal value at $r=b(s)$
for fixed $s$. 
Because of the restriction (\ref{<2}), we have
%the condition (\ref{b(s)}) uniquely determines $b(s)$.
\begin{equation}\label{b(s)}
[tb(s)]=s\hbox{ or }s-1.
\end{equation}
Here $[x]$ is the integer part of $x$.
We construct vectors in the spaces $M_{r,s}$ ($1\leq r\leq p-1$)
by applying the Fourier components $\phi^\pm_n$ of the $(2,1)$-primary field
(\ref{phi21}) to $|b(s),s\rangle$.

To each triple of integers $(r,r',r'')$ satisfying the conditions
$1\leq r,r',r''\leq p-1$, $r'=r\pm1$ and $r''=r'\pm1$ we associate 
a local weight $w(r,r',r'')$:
\begin{eqnarray}
&&w(r,r\pm1,r\pm2)=\frac t2,\\
&&w(r,r+1,r)=2-\frac t2+[tr]-tr,\\
&&w(r,r-1,r)=1-\frac t2-[tr]+tr.
\end{eqnarray}
We have
\begin{eqnarray}
&&w(r,r',r'')\geq0,\\
&&w(r,r',r'')\equiv\Delta_{r'',s}+\Delta_{r,s}-2\Delta_{r',s}\bmod\Z,\\
&&w(r,r',r'')=w(p-r,p-r',p-r'').
\end{eqnarray}

A sequence of integers ${\bf r}=(r_0,r_1,\ldots,r_L)$ satisfying the conditions
$1\leq r_i\leq p-1$ and $r_{i+1}=r_i\pm1$ is called a one-dimensional
configuration of length $L$. We denote by $C^{(L)}_{a,c}$ the set of
one-dimensional configurations of length $L$ satisfying $r_0=a$ and $r_L=c$.
By the definition $C^{(L)}_{a,c}$ is an empty set unless $L\equiv a-c\bmod2$.

Let $L\geq0$ be an integer. We define type $s$ admissible monomials of length
$L$, and associated one-dimensional configurations of length $L$.

A type $s$ monomial of length $L$ is a sequence
$m=(\sigma_1,m_1;\ldots,;\sigma_L,m_L)$ of signs $\sigma_i=\pm$
(or $\sigma_i=\pm1)$ and rational numbers $m_i\in{\bf Q}$ ($1\leq i\leq L$).
The associated one-dimensional configuration ${\bf r}=r(m)$
of length $L$ is defined by
\begin{equation}
r(m)_L=b(s)\hbox{ and }r(m)_{i-1}=r(m)_i+\sigma_i\quad(2\leq i\leq L).
\end{equation}

A void sequence is a type $s$ monomial of length $0$. We denote it by
$\emptyset$. The monomial $\emptyset$ is the unique type $s$ monomial of length
$0$. It is an admissible monomial by definition. The associated one-dimensional
configuration $r(\emptyset)$ is such that $r(\emptyset)_0=b(s)$.

If $L\geq 1$, a type $s$ monomial $m$ is admissible if and only if
the following conditions, where $1\leq i\leq L-1$, are valid.
\begin{eqnarray}
&&1\leq r(m)_i\leq p-1,\label{R1}\\
&&-m_i\in\Delta_{r(m)_{i-1},s}-\Delta_{r(m)_i,s}+\Z,\label{R0}\\
&&-m_L\in\Delta_{r(m)_{L-1},s}-\Delta_{b(s),s}+\Z_{\geq0},\label{R2}\\
&&-m_i+m_{i+1}\in w(r(m)_{i-1},r(m)_i,r(m)_{i+1})+\Z_{\geq0}.\label{R3}
\end{eqnarray}

With each admissible monomial
$m=(\sigma_1,m_1;\ldots,;\sigma_L,m_L)$ of type $s$
we associate the product of the Fourier components:
\begin{equation}\label{MONO}
\Phi(m)=\phi^{\sigma_1}_{m_1}\ldots\phi^{\sigma_L}_{m_L},
\end{equation}
and define a vector $|m\rangle\in M_{r,s}$ where $r=r(m)_0$ by
\begin{equation}
|m\rangle=\Phi(m)|b(s),s\rangle.
\end{equation}
The degree of the vector $|m\rangle$ is given by
\begin{equation}\label{d(m)}
d(m)=\Delta_{b(s),s}-\sum_{i=1}^Lm_i,
\end{equation}
This is consistent with the degree $d$ in $M_{r,s}=\oplus_d(M_{r,s})_d$ given
by the operator $L_0$:
\begin{equation}
(M_{r,s})_d=\{v\in M_{r,s}|L_0v=dv\}.
\end{equation}
The operator $\phi^\sigma_n$ changes the degree by $-n$.

Starting from the vector $|b(s),s\rangle$, we create vectors
successively by the operators $\phi^{\sigma_i}_
{m_i}$.
The conditions (\ref{R1}) and (\ref{R0}) ensure that the operator
$\phi^{\sigma_i}_{m_i}$ is acting as
$M_{r(m)_i,s}\rightarrow M_{r(m)_{i-1},s}$, the condition (\ref{R2})
must hold if the vector $\phi^{\sigma_L}_{m_L}|b(s),s\rangle$
is non-zero, and the condition (\ref{R3}) requires that
the increment from $i+1$ to $i$ of the degree differences caused
by the operators $\phi^{\sigma_{i+1}}_{m_{i+1}}$ and $\phi^{\sigma_i}_{m_i}$
is at least $w(r(m)_{i-1},r(m)_i,r(m)_{i+1})$.

We denote by $B_{r,s}$ the set of type $s$ admissible monomials
$m$ such that $r(m)_0=r$. In $B_{r,s}$ the length of elements varies.
We will prove
\begin{prop}\label{prop:SPANNING}
The set of vectors $\{|m\rangle|m\in B_{r,s}\}$ is a spanning set
of vectors in $M_{r,s}$.
\end{prop}

We prepare some notation before starting the proof of Proposition
\ref{prop:SPANNING}.

We denote by $\tilde B_{r,s}$ the set of type $s$ monomials satisfying
(\ref{R1}), (\ref{R0}), (\ref{R2}) and $r(m)_0=r$. We drop the conditions
(\ref{R3}) from $B_{r,s}$. The vectors $|m\rangle$ are
defined for monomials $m\in\tilde B_{r,s}$ as well. We have
\begin{prop}\label{VOA}
The set of vectors $\{|m\rangle|m\in\tilde B_{r,s}\}$ is a spanning set
of vectors in $M_{r,s}$.
\end{prop}
%We will prove
%\begin{prop}\label{SPANNING}
%The set of vectors %$\{|m\rangle|m\in B_{r,s},1\le r\le p-1\}$ 
%$\phi^{\sigma_1}_{m_1}\ldots\phi^{\sigma_L}_{m_L}\ket{b(s),s}$,
%where $L\ge 0$ and $\sigma_i,m_i$ run over all possible values, 
%is a spanning set of vectors in $\oplus_{r=1}^{p-1}M_{r,s}$.
%\end{prop}
\begin{proof}
Let $V\subset \oplus_{1\le r\le p-1}M_{r,s}$ be the linear 
span of the vectors 
$\phi^{\sigma_1}_{m_1}\ldots\phi^{\sigma_L}_{m_L}\ket{b(s),s}$.
The proposition will follow if we show
\bea
&&\ket{r,s}\in V \qquad (1\le r\le p-1),
\label{chk1}
\\
&&L_nV\subset V\qquad (n\in\Z).
\label{chk2}
\ena
Assertion \eqref{chk1} is clear from 
\eqref{normalization_of_primary}. 
Assertion \eqref{chk2} for $n\ge 0$ follows from 
the intertwining relation \eqref{intertwiner} 
and the highest condition $L_n\ket{r,s}=0$ ($n\ge 0$). 
%By the same argument, t
To verify \eqref{chk2} for $n< 0$,  
it is enough to show that $L_n\ket{b(s),s}\in V$ for $n=-1,-2$.  
This can be seen by the following formula
obtained by using \eqref{normalization_of_primary},
\eqref{intertwiner}:
\be
&&L_{-1}\ket{r,s}=\pm\frac{(r\mp1)t-s\pm1}{t}
\phi^\pm_{\Delta_{r\mp1,s}-\Delta_{r,s}-1}\ket{r\mp1,s},
\\
&&
\bigl(L_{-2}\pm\frac{t}{(r\mp1)t-s\pm1}L_{-1}^2\bigr)
\ket{r,s}
\\
&&
=\pm\frac{2}{t}\bigl((r\mp1)t-s\pm2\bigr)
\phi^\pm_{\Delta_{r\mp1,s}-\Delta_{r,s}-2}\ket{r\mp1,s}.
\\
\en
The proof is over.
\end{proof}
%\end{document}
%\begin{proof}
%Let ${\mathcal V}_{2,1}$ be the algebra generated by the Fourier components of
%the primary field $\phi_{2,1}(z)$.
%Because of
%(\ref{normalization_of_primary})
%we obtain all the vectors $|r,s\rangle$ ($1\leq r\leq p-1$)
%in the ${\mathcal V}_{2,1}$ orbit of $|b(s),s\rangle$.
%Since ${\mathcal V}_{2,1}$ contains the Virasoro algebra, we obtain all the
%spaces $M_{r,s}$ ($1\leq r\leq p-1$) by the action of ${\mathcal V}_{2,1}$.
%Namely, we obtain any vector in $\oplus_{r=1}^{p-1}M_{r,s}$
%as a linear combination of the vectors $|m\rangle$ where the monomials $m$
%satisfy the conditions (\ref{R1}) and (\ref{R0}).
%If a monomial $m$ breaks the condition (\ref{R2}),
%then we have $|m\rangle=0$. Therefore, we can discard it. Thus,
%we obtain the spanning set $\tilde B_{r,s}$.
%\end{proof}

The vectors in Proposition \ref{VOA} are not linearly independent.
We want to discard those which are linearly dependent on others.
In order to do this procedure systematically, we define a partial ordering of
the set $\tilde B_{r,s}$. We will show that if a monomial $m\in\tilde B_{r,s}$
is non-admissible the vector $|m\rangle$ can be written as a linear
combination of the vectors associated with monomials smaller than $m$.

Let $m,m'\in\tilde B_{r,s}$, and let $L,L'$ be the lengths of $m,m'$,
respectively. We write $m<m'$ if and only if $L<L'$, or
$L=L'$ and there exists $1\leq l\leq L$ such that
$m_i=m'_i$ for all $1\leq i\leq l-1$ and $m_l<m'_l$.
Note that the above ordering is not a total order because we
do not compare $m,m'$ if they have the same length $L$ and $m_i=m'_i$
for all $1\leq i\leq L$, but $r(m)\not=r(m')$.
Since the number of monomials of a fixed degree is finite, it follows that the
vectors associated with non-admissible monomials are contained in the
linear span of those associated with admissible monomials.

\begin{prop}\label{SPAN}
If a monomial $m\in\tilde B_{r,s}$ does not satisfy $(\ref{R3})$,
then it can be written as a linear combination of vectors 
corresponding to smaller monomials in $\tilde B_{r,s}$.
\end{prop}
Before starting the proof we prepare some notation and a technical lemma.

Recall that $\rho=rt-s$, $\lambda_\pm=t/4\pm\rho/2$. We set
\begin{eqnarray}
&&\delta=
\begin{cases}
0\hbox{ if $n_1+n_2$ is even;}\\
1\hbox{ if $n_1+n_2$ is odd,}
\end{cases}\\
&&w_+=w(r,r+1,r)=2-t/2+[\rho]-\rho,\\
&&w_-=w(r,r-1,r)=1-t/2-[\rho]+\rho.
\end{eqnarray}
Note that $1<w_-+w_+<2$.

We expand the hypergeometric functions which appear in (\ref{1.27}):
\begin{eqnarray}
f^\pm(z)&=&(1-z)^{-2+\frac{3t}2}F(-1+t,-1+t\mp\rho,1\mp\rho;z)\label{f}\\
&=&(1-z)^{1-\frac t2}F(2-t,2-t\mp\rho,1\mp\rho;z)\nonumber\\
&=&\sum_{n\in\Z}f^\pm_nz^n.\nonumber
\end{eqnarray}
For (\ref{1.28}), we have
\begin{eqnarray}
g^\pm(z)&=&(z\partial-\lambda_\pm)\left((1-z)f^\pm(z)\right)\label{g}\\
&=&\sum_{n\in\Z}g^\pm_nz^n.\nonumber
\end{eqnarray}
Note that if $n<0$, we have $f^\pm_n=g^\pm_n=0$.

\begin{lem}\label{MATLEM}
Define a $(2n+\delta)\times(2n+\delta)$ matrix
\begin{equation}
M_{n,\delta}(\rho)_{j,k}=
\begin{cases}
f^+_{j+k-n-1-\delta}+f^+_{-j+k-n}\hbox{ if $1\leq j\leq n+\delta$;}\\
g^+_{j+k-2n-1-\delta}+g^+_{-j+k+1}\hbox{ if $n+\delta+1\leq j\leq2n+\delta$.}
\end{cases}
\end{equation}
The matrix $M_{n,\delta}(\rho)$ is non-singular.
\end{lem}
\begin{proof}
Set
\begin{equation}
D_{n,\delta}(\rho)={\rm det}\,M_{n,\delta}(\rho).
\end{equation}
We will show that up to a non-zero constant the determinant factorizes:
\begin{equation}\label{DET}
D_{n,\delta}(\rho)=\hbox{const.}\frac{\prod_{i=1}^{2n-1+\delta}
\{(\rho+t-i-1)(\rho-t-i+2)\}^{n-\left[\frac{i-\delta}2\right]}}
{\prod_{i=1}^{2n-1+\delta}(\rho-i)^{2n-i+\delta}}
\end{equation}
Since the determinant never vanishes when $2\leq r\leq p-2$, the proof will
be thus finished.

Set
\begin{equation}
f^+_n=\frac{a_n(\rho)}{\prod_{i=1}^n(\rho-i)},\quad
g^+_n=\frac{b_n(\rho)}{\prod_{i=1}^n(\rho-i)}.
\end{equation}
{}From (\ref{f}) and (\ref{g}) we see that $a_n(\rho)$ and $b_n(\rho)$
are polynomials in $\rho$ of degree $n$ and $n+1$, respectively. This implies
that if we multiply $D_{n,\delta}(\rho)$ by the denominator in the right
hand side of (\ref{DET}), we have a polynomial in $\rho$ of degree
$2n(n+\delta)$. Therefore, in order to prove (\ref{DET}) it is enough to show
that for $1\leq i\leq 2n-1+\delta$ by specializing at $\rho=-t+i+1$ or
$\rho=t+i-2$, the corank of the matrix $M_{n,\delta}(\rho)$ becomes at least
$n-[(i-\delta)/2]$.

We prove the statement above for the case $\delta=0$ and $\rho=-t+i+1$. Three
other cases are similar. From (\ref{f}) we see that at the special values
$\rho=-t+i+1$ ($i\geq1$) the series $F(z)=(1-z)^{-1+\frac t2}f^+(z)$ becomes
a reciprocal polynomial of degree $i-1$. Namely, we have
$F(z)=\sum_{j=0}^{i-1}F_jz^j$ where $F_j=F_{i-1-j}$. Suppose that $P(z)$ is
a reciprocal polynomial of degree $l$. Then, the polynomials $(1+z)P(z)$ and
$\left((1-z)z\partial+lz\right)P(z)$ are also reciprocal and their degrees are
$l+1$ and $l$. From this remark and (\ref{g}) it follows that the series
$G(z)=(1-z)^{-1+\frac t2}g^+(z)$ is a reciprocal polynomial of degree $i$.

Let $(1-z)^{-1+\frac t2}=\sum_{j\in\Z}c_jz^j$ be the expansion at $z=0$.
If we replace $f^+_j,g^+_j$ by the coefficients $F_j,G_j$ of $F(z),G(z)$,
the matrix $M_{n,\delta}(\rho)$ is left-multiplied by $C=(C_{jk})$
where $C_{jk}=c_{k-j}$. Therefore, it is enough to prove the statement when
the matrix $M_{n,\delta}(\rho)$ is given by $F_j,G_j$ instead of $f^+_j,g^+_j$.
On closer inspection, the matrix thus obtained proves to have the symmetries
\begin{equation}
(CM_{n,0}(-t+i+1))_{j,k}=(CM_{n,0}(-t+i+1))_{j,2n+i-k}
\end{equation}
for $1\leq j,k\leq n$ and $1\leq2n+i-k\leq n$. The assertion
follows from this.

To show that the constant in \eqref{DET} is non-zero, we specialize
$\rho=-1$ and $t=2$ which produces a simple matrix with determinant
$\pm 2$.
\end{proof}

{\it Proof of Proposition \ref{SPAN}.}
We start from Proposition \ref{VOA}.
Suppose that a monomial $m\in\tilde B_{r,s}$ does not satisfy (\ref{R3}).
Then, we have $L\geq2$, and there exists $1\leq i\leq L-1$ such that
\begin{equation}
-m_i+m_{i+1}<w(r(m)_{i-1},r(m)_i,r(m)_{i+1}).
\end{equation}
For the proof of Proposition \ref{SPAN} it is enough to show that the vector 
$\phi^{\sigma_i}_{m_i}\phi^{\sigma_{i+1}}_{m_{i+1}}|v\rangle$ where
$|v\rangle=\phi^{\sigma_{i+2}}_{m_{i+2}}\cdots\phi^L_{m_L}|b(s),s\rangle$
can be written as a linear combination of vectors of the form
$\phi^{\sigma'}_{n'}\phi^{\sigma''}_{n''}|v\rangle$ satisfying
$\sigma_i+\sigma_{i+1}=\sigma'+\sigma''$,
$m_i+m_{i+1}=n'+n''$ and $n'<m_i$.

We use the quadratic relations given in Proposition \ref{QUADREL} to rewrite
the product $\phi^{\sigma_i}_{m_i}\phi^{\sigma_{i+1}}_{m_{i+1}}$.
There are six cases:

{\it Case 1}\,: $(r(m)_{i-1},r(m)_i,r(m)_{i+1})=(r,r+1,r+2)$

{\it Case 2}\,: $(r(m)_{i-1},r(m)_i,r(m)_{i+1})=(r,r-1,r-2)$

{\it Case 3}\,: $(r(m)_{i-1},r(m)_i,r(m)_{i+1})=(1,2,1)$

{\it Case 4}\,: $(r(m)_{i-1},r(m)_i,r(m)_{i+1})=(p-1,p-2,p-1)$

{\it Case 5}\,: $(r(m)_{i-1},r(m)_i,r(m)_{i+1})=(r,r+1,r)$ ($2\leq r\leq p-2$)

{\it Case 6}\,: $(r(m)_{i-1},r(m)_i,r(m)_{i+1})=(r,r-1,r)$ ($2\leq r\leq p-2$)

By the symmetry $(r,s)\leftrightarrow(p-r,p'-s)$,
Case 2 is equivalent to Case 1, and Case 4 is equivalent to Case 3.

Case 1 is discussed in Introduction. We will not repeat the argument.
%The quadratic relation of the form (\ref{Q3})
%has appeared \cite{LWP} in the representation theory of the affine $sl_2$
%modules. The argument in \cite{} applies to the proof of Case 1 without
%alteration.
We can ignore the right hand side of the quadratic identities (\ref{1.27}) and
(\ref{1.28}), because it contributes with only smaller terms to the corresponding
monomials. If we forget the right hand side, the quadratic relation
(\ref{1.27})
for $r=1$ (see (\ref{1.18})) is similar to (\ref{1.29}) (see (\ref{1.20})).
The difference is in the sign of the second term, and that the power
$(z_1-z_2)^{-t/2}$ in (\ref{1.20}) is replaced  by $(z_1-z_2)^{-2+3t/2}$ in
(\ref{1.18}). These powers are related to the weight $w$ in (\ref{R3}).
The effect of the change of sign is that $w(1,2,1)=3-3t/2$ while
$w(r,r+1,r+2)=t/2$.

Cases 5 and 6 are combined in the relations (\ref{1.27}) and (\ref{1.28}),
and thus the proof is more involved.

By the symmetry $(r,s)\leftrightarrow(p-r,p'-s)$,
without loss of generality we can assume that
\begin{equation}
\rho=rt-s>0.
\end{equation}
We have the Fourier series expansions
\begin{eqnarray}
&&\phi^{(r,s,r\pm1,s)}_{2,1}(z)
=\sum_{n\in\Z+\Delta_{r\pm1,s}-\Delta_{r,s}}
\phi^{(r,r\pm1)}_nz^{-n-\Delta_{2,1}},\label{+}\\
&&\phi^{(r\pm1,s,r,s)}_{2,1}(z)
=\sum_{n\in\Z-\Delta_{r\pm1,s}+\Delta_{r,s}}
\phi^{(r\pm1,r)}_nz^{-n-\Delta_{2,1}}.\label{-}
\end{eqnarray}
Here we fix $s$ and suppress it from the notation.
Substituting (\ref{+}) and (\ref{-}) in the left hand side of (\ref{1.27}) and
taking the coefficient of each monomial in $z_1,z_2$, we obtain an infinite 
linear combination of the operator products
$\phi^{(r,r\pm1)}_{n_1}\phi^{(r\pm1,r)}_{n_2}$. Modulo a constant in the right
hand side, we obtain a relation among the monomials 
$\phi^{(r,r\pm1)}_{n_1}\phi^{(r\pm1,r)}_{n_2}$. Similarly, we obtain another set
of relations from (\ref{1.28}). Each of these relations contains infinitely
many terms of the form $\phi^{(r,r\pm1)}_{n_1}\phi^{(r\pm1,r)}_{n_2}$.
The sum $n_1+n_2$ is constant for all the terms which appear in one relation,
and moreover, $n_2$ is bounded from below. Therefore, acting on each vector
in $M_{r,s}$, only finitely many terms create non-zero vectors.
Namely, the infinite quadratic relation for the operators gives a
finite linear relation among the vectors created by them.

Let us write explicitly the relation. We fix $n_1+n_2$. In order to
simplify the notation, we write
\begin{equation}
\Phi^\pm_{n_2-n_1}=\phi^{(r,r\pm1)}_{n_1}\phi^{(r\pm1,r)}_{n_2}.
\end{equation}
since the sum $n_1+n_2$ is fixed, there is no ambiguity when we write
only the difference $n_2-n_1$. Thus, $\Phi^\pm_m$ is defined for
$m\equiv n_1+n_2\bmod2$.

For $a\in\Z$ we set
\begin{eqnarray}
A_a=\sum_\pm\sum_{i\in\Z}(f^\pm_i+f^\pm_{i-2a-1+\delta})
\Phi^\pm_{-2a-2\lambda_\pm+2i+\delta},\label{A_a}\\
B_a=\sum_\pm\sum_{i\in\Z}(g^\pm_i+g^\pm_{i-2a+\delta})
\Phi^\pm_{-2a-2\lambda_\pm+2i+\delta}.\label{B_a}
\end{eqnarray}

{}From (\ref{1.27}) we obtain the relations $A_a\equiv0$ ($a\in\Z_{\geq0}$),
and 
from (\ref{1.28}) $B_a\equiv0$ ($a\in\Z_{\geq\delta}$).
Here $\equiv$ means that we ignore constant terms. In the proof below we always
ignore constant terms.

In accordance with admissibility of monomials, we say $\Phi^\pm_n$ is
admissible if and only if
\begin{equation}
n\in w_\pm+\Z_{\geq0}.
\end{equation}
We say also $\Phi^\sigma_n$ is larger than $\Phi^{\sigma'}_{n'}$
if and only if $n<n'$ (not opposite).

Our goal is to show that if $\Phi^\sigma_n$ is non-admissible then it can
be written as an infinite linear combination of $\Phi^\pm_{n'}$
such that $n'\geq w_\pm$. Since $w_\pm-1<w_\mp$, it implies that the product
$\phi^{\sigma_i}_{m_i}\phi^{\sigma_{i+1}}_{m_{i+1}}$ can be replaced
by a linear combination of $\phi^{\sigma'}_{n'}\phi^{\sigma''}_{n''}$
such that $n'<m_i$ (modulo a constant).

Let us prove the above statement.

We set
\begin{equation}
n_\pm(a)=-2a-2\lambda_\pm+\delta.
\end{equation}
The largest term among those $\Phi^\pm_n$ which appear in $A_a$
($a\in\Z_{\geq0}$) or $B_a$ ($a\in\Z_{\geq\delta}$)
is $\Phi^\pm_{n_\pm(a)}$:
\begin{eqnarray}
A_a&=&f^+_0\Phi^+_{n_+(a)}+\cdots\\
&&+f^-_0\Phi^-_{n_-(a)}+\cdots,\nonumber\\
B_a&=&g^+_0\Phi^+_{n_+(a)}+\cdots\\
&&+g^-_0\Phi^-_{n_-(a)}+\cdots.\nonumber
\end{eqnarray}
Since the matrix
\begin{equation}
\left(
\begin{matrix}
f^+_0&f^-_0\\
g^+_0&g^-_0\\
\end{matrix}
\right)=
\begin{pmatrix}1 &1 \\
-\lambda^+&-\lambda_-\\
\end{pmatrix}
\end{equation}
is non-degenerate, we can replace $\Phi^\pm_{n_\pm(a)}$ by smaller terms
$\Phi^\sigma_n$, i.e., $n>n_\sigma(a)$. We do replace them while
both $\Phi^+_{n_+(a)}$ and $\Phi^-_{n_-(a)}$ are non-admissible, i.e.,
$n_\pm(a)<w_\pm$. This is equivalent to $a\geq N$ where
\begin{equation}
N=
\begin{cases}
\frac{[\rho]}2+\delta\hbox{ if $[\rho]$ is even;}\\
\frac{[\rho]+1}2\hbox{ if $[\rho]$ is odd.}
\end{cases}
\end{equation}
For this value, we have
\begin{eqnarray}
n_+(N)&=&
\begin{cases}
w_+-2[\rho]-2-\delta\hbox{ if $[\rho]$ is even;}\\
w_+-2[\rho]-3+\delta\hbox{ if $[\rho]$ is odd,}
\end{cases}
\\
n_-(N)&=&
\begin{cases}
w_--1-\delta\hbox{ if $[\rho]$ is even;}\\
w_--2+\delta\hbox{ if $[\rho]$ is odd.}
\end{cases}
\end{eqnarray}
Note, in particular, that $n_-(N)<w_-$ and $n_-(N)+2\ge w_-$.
Therefore, the non-admissible terms $\Phi^\sigma_n$ with $\sigma=-$,
and those with $\sigma=+$ and $n\in n_+(N)+2\Z_{\leq0}$ can be
replaced by the terms $\Phi^+_n$ ($n\in n_+(N)+2\Z_{>0}$) and the
admissible terms $\Phi^-_n$ ($n\geq w_-$).

The remaining  non-admissible terms are
\begin{eqnarray}
&&\Phi^+_{w_+-2a+\delta}\quad(1\leq a\leq 2N-\delta)
\hbox{ if $[\rho]$ is even;}\\
&&\Phi^+_{w_+-2a+1-\delta}\quad(1\leq a\leq 2N-\delta)
\hbox{ if $[\rho]$ is odd.}
\end{eqnarray}

We want to eliminate these terms by
using the relations $A_a\equiv0$ ($0\leq a\leq N-1$) and
$B_a\equiv0$ ($\delta\leq a\leq N-1$). For this we need to show that the
$(2N-\delta)\times(2N-\delta)$ matrix whose elements are the coefficients of
the non-admissible terms in the these relations, is non-degenerate.
If $[\rho]=2n$ ($n\geq1$), we put the coefficients of
$\Phi^+_{w_+-4N+3\delta},\ldots,\Phi^+_{w_+-2+\delta}$ in the first,..., the
$(2N-\delta)$-th column of the matrix, respectively.
If $[\rho]=2n-1+2\delta$ ($n\geq1-\delta$), we put the coefficients of
$\Phi^+_{w_+-4N+1+\delta},\ldots,\Phi^+_{w_+-1-\delta}$ in the first,...,
the $(2N-\delta)$-th column of the matrix, respectively.
We put the coefficients in
$A_0,\ldots,A_{N-1},B_\delta,\ldots,B_{N-1}$ in the first,..., the
$(2N-\delta)$-th row of the matrix, respectively.
If $\delta$ is fixed, two matrices corresponding to $[\rho]=2n,2n-1+2\delta$
are the same. We denote it by $M_{n,\delta}(\rho)$. The size of
$M_{n,\delta}(\rho)$ is $2n-\delta$. This matrix is nothing but
$M_{n,\delta}(\rho)$ given by Lemma \ref{MATLEM}. The proof is over. \qed

\subsection{One-dimensional configuration sums}
In this section we calculate the character of the spanning set
$B_{r,s}$:
\begin{equation}
\tilde\chi_{r,s}(q)=\sum_{m\in B_{r,s}}q^{d(m)},
\end{equation}
where $d(m)$ is given by (\ref{d(m)}).
The character $\tilde\chi_{r,s}(q)$ is expressed as follows by using
one-dimensional configuration sums:
\begin{prop}
\begin{eqnarray}
&&\tilde\chi_{r,s}(q)=q^{\Delta_{b(s),s}}\Bigl(\delta_{r,b(s)}
+\sum_{L\geq1}\frac1{(q)_L}
\sum_{{\bf r}\in C^{(L)}_{r,b(s)}}q^{d_{\bf r}}\Bigr),\\
&&d_{\bf r}=
L(\Delta_{r_{L-1},s}-\Delta_{b(s),s})+\sum_{i=1}^{L-1}iw(r_{i-1},r_i,r_{i+1}).
\end{eqnarray}
\end{prop}
\begin{proof}
The degree of the vector $|b(s),s)\rangle$ is $\Delta_{b(s),s}$.
The operator product $\Phi(m)$ (see (\ref{MONO}) adds some degree to it.
The minimal degree added by those $\Phi(m)$ which have
the same associated one-dimensional configuration, i.e.,
$r(m)=\bf r$, is equal to $d_{\bf r}$. The fluctuation from the minimum is
added up to the factor $1/(q)_L$.
\end{proof}
Let $a,b,c$ be integers satisfying
\begin{equation}
1\leq a,b,c\leq p-1,b=c\pm1.
\end{equation}
We introduce the one-dimensional configuration sum
\begin{equation}
Y^{(L)}_{a,b,c}(q)=\sum_{{\bf r}\in C^{(L+1)}_{a,c},r_L=b}
q^{\sum_{i=1}^Liw(r_{i-1},r_i,r_{i+1})}.
\end{equation}
Note that
\begin{equation}
Y^{(0)}_{a,b,c}(q)=\delta_{a,b}.
\end{equation}
We have
\begin{equation}\label{1DSUM}
\tilde\chi_{r,s}(q)=\delta_{r,b(s)}q^{\Delta_{b(s),s}}
+\sum_{L\geq1}\frac{q^{\Delta_{b(s),s}}}{(q)_L}
\sum_{\tau=\pm1} q^{L(\Delta_{b(s)+\tau,s}-\Delta_{b(s),s})}
Y^{(L-1)}_{r,b(s)+\tau,b(s)}.
\end{equation}
Similar sums appeared in the calculation of local height probabilities
of the eight vertex solid-on-solid model (see \cite{ABF,BF}).
In fact, there is a connection (first observed in \cite{Huse})
between one-dimensional configuration sums and characters in
conformal field theories. The Virasoro character $\chi_{r,s}$ in
the minimal series for general $(p,p')$ (i.e., without the restriction $p'<2p$)
is obtained (see (6) in \cite{FLPW}) in the limit $L\rightarrow\infty$ of
the one-dimensional configuration sum where the local weights are given by
\begin{eqnarray}\label{BF}
&&\tilde w(r,r\pm1,r\pm2)=\frac12,\\
&&\tilde w(r,r+1,r)=-[r(p'-p)/p'],\\
&&\tilde w(r,r-1,r)=[r(p'-p)/p'].
\end{eqnarray}
Set
\begin{equation}
g(r)=-\frac14(t-1)r^2+\frac14r.
\end{equation}
If we modify our weight $w$ to $w'$ by the gauge transformation
\begin{equation}\label{GAUGE}
w'(r,r',r'')=w(r,r',r'')-1+g(r)-2g(r')+g(r''),
\end{equation}
we obtain
\begin{eqnarray}\label{OURS}
&&w'(r,r\pm1,r\pm2)=-\frac12,\\
&&w'(r,r+1,r)=[r(p'-p)/p],\\
&&w'(r,r-1,r)=-[r(p'-p)/p].
\end{eqnarray}
Since
\[
\sum_{i=1}^Liw(r_{i-1},r_i,r_{i+1})=\sum_{i=1}^Liw'(r_{i-1},r_i,r_{i+1})
+\frac{L(L+1)}2-g(r_0)+(L+1)g(r_L)-Lg(r_{L+1}),
\]
the gauge transformation does not essentially change the one-dimensional
configuration sum. The expressions (\ref{BF}) and (\ref{OURS})
are very similar if we change the sign. However, they are not equal
for the same values of $p,p'$. In fact, the way of connecting
the one-dimensional configuration sums for $\tilde w$ to the
Virasoro characters is very different from our way of connecting those
for $w$ to the same Virasoro characters. The former uses the limit of $L$,
while in our formula the parameter $L$ appears as a summation variable.

By a routine calculation (see, e.g., \cite{ABF}), we have the following result.
\begin{prop}
For integers $a,b,c$ satisfying $b=c\pm1$, we define
\begin{equation}
X^{(L)}_{a,b,c}(q)=\sum_{
{{\bf r}=(r_0,\ldots,r_L,r_{L+1})\atop
r_{i+1}=r_i\pm1(0\leq i \leq L)}
\atop r_0=a,r_L=b,r_{L+1}=c}
q^{\sum_{i=1}^Liw(r_{i-1},r_i,r_{i+1})}.
\end{equation}
Then, we have
\begin{equation}
X^{(L)}_{a,b,b\pm1}(q)=q^C\left[L\atop\frac{L\pm(a-b)}2\right]
\end{equation}
where
\begin{eqnarray}
&&C=\frac14\Bigl\{L(L+1)\mp L+(t-1)(a^2-b^2)-a+b\nonumber\\
&&+L(t-1)(\pm2b+1)+(1\mp2[(t-1)c])(L\pm(a-b))\Bigr\},\nonumber
\end{eqnarray}
and
\begin{equation}\label{YX}
Y^{(L)}_{a,b,b\pm1}(q)=\sum_\varepsilon\varepsilon\sum_{n\in\Z}
X^{(L)}_{\varepsilon a+2np,b,b\pm1}(q).
\end{equation}
\end{prop}
We omit the proof of this proposition.
We only remark that the gauge transformation (\ref{GAUGE}) makes
the computation shorter, and that
the summation (\ref{YX}) realizes the boundary conditions
\begin{equation}
Y^{(L)}_{a,0,1}(q)=Y^{(L)}_{a,p,p-1}(q)=0\hbox{ for }1\leq a\leq p-1.
\end{equation}

\subsection{Fermionic formulas of Virasoro characters}
We denote the corresponding character by $\chi_{r,s}$.
It is given by the following bosonic formula.
\begin{equation}\label{CHAR}
\chi_{r,s}(q)=\frac{q^{\Delta_{r,s}}}{(q)_\infty}
\left(\sum_{n\in\Z}q^{pp'n^2+(p'r-ps)n}-
\sum_{n\in\Z}q^{pp'n^2+(p'r+ps)n+rs}\right).
\end{equation}
We rewrite it by using
\begin{lem}
For all integer $l$, we have
\begin{equation}
\frac1{(q)_\infty}=\sum_{m\in\Z}\frac{q^{m^2-ml}}{(q)_m(q)_{m+l}}.
\end{equation}
\end{lem}
Here, we set
\begin{equation}
\frac1{(q)_m}=0\hbox{ if }m<0.
\end{equation}

Applying the lemma to the first/second sum in (\ref{CHAR}) by setting
\begin{equation}
l=\mp r+b+2pn,
\end{equation}
and changing the summation over $m\in\Z$ to the summation over
\begin{equation}
L=2m+2pn+l\in\mp r+b +2\Z,
\end{equation}
we obtain
\begin{prop}The Virasoro character $\chi_{r,s}$ can be written as
follows where $b$ is any integer.
\begin{eqnarray}
&&\chi_{r,s}(q)=q^{\Delta_{r,s}}
\sum_{L\geq0}\frac{q^{\frac{L^2-(r-b)^2}4}}{(q)_{\scriptscriptstyle L}}
\sum_{n\in\Z}\Biggl(q^A
\left[\scriptstyle L\atop\scriptstyle\frac{L-r+b}2-pn\right]
-q^B\left[\scriptstyle L\atop\scriptstyle\frac{L-r-b}2-pn\right]\Biggr)
\label{FERM}\\
&&\hbox{where}\nonumber\\
&&A=p(p'-p)n^2+\{(p'-p)r-p(s-b)\}n,\nonumber\\
&&B=p(p'-p)n^2+\{(p'-p)r+p(s-b)\}n+r(s-b).\nonumber
\end{eqnarray}
\end{prop}
We will identify the formula (\ref{FERM}) where $b=b(s)$ with (\ref{1DSUM}).
As a result, we see that if $b=b(s)$, for each $L$
the sum over $n$ in (\ref{FERM}) is a series with non-negative coefficients.
\begin{prop}
We have the equality
\begin{equation}
\tilde\chi_{r,s}(q)=\chi_{r,s}(q).
\end{equation}
\end{prop}
\begin{proof}
We set $b=b(s)$ in (\ref{FERM}). There are two cases:
Case\,(i) $[tb(s)]=s$; Case\,(ii) $[tb(s)]=s-1$.
In Case\,(i), we rewrite (\ref{FERM}) by using the identities
\begin{eqnarray}
\left[L\atop\frac{L-r+b}2-pn\right]=
q^{\frac{L-r+b}2-pn}\left[L-1\atop\frac{L-r+b}2-pn\right]+
\left[L-1\atop\frac{L+r-b}2+pn\right],\\
\left[L\atop\frac{L-r-b}2-pn\right]=
\left[L-1\atop\frac{L-r-b}2-pn\right]+
q^{\frac{L+r+b}2+pn}\left[L-1\atop\frac{L+r+b}2+pn\right].
\end{eqnarray}
For each $L\geq1$, we thus obtain four terms in (\ref{FERM}).
One can check these four terms are equal to the terms in (\ref{1DSUM})
and (\ref{YX})
corresponding to $(\varepsilon,\tau)=(1,1),(1,-1),(-1,-1),(-1,1)$,
respectively. The proof for Case\,(ii) is similar.
\end{proof}
As a corollary, we have
\begin{thm}
The set of vectors $|m\rangle$ $(m\in B_{r,s})$ is a basis of
the representation space $M_{r,s}$.
\end{thm}
\section{The $\slth$-case}

In this section we construct a basis of the 
integrable irreducible highest weight modules of the affine Lie algebra $\slth$. 
The construction is quite similar to the case of Virasoro modules. 
In the $\slth$-case we use the vertex operators (see \eqref{eq:def-VO} below) 
instead of the $(2, 1)$-primary fields. 

\subsection{Preliminaries}  
First we introduce some notation on the affine Lie algebra $\slth$ and 
its representations. 

Denote by $E, F$ and $H$ the generators of $\mathfrak{sl}_{2}$ given by 
\be 
E=\left( 
\begin{array}{cc} 
0 & 1 \\ 0 & 0 
\end{array} 
\right), \quad 
F=\left( 
\begin{array}{cc} 
0 & 0 \\ 1 & 0 
\end{array} 
\right), \quad 
H=\left( 
\begin{array}{cc} 
1 & 0 \\ 0 & -1 
\end{array} 
\right). 
\en 
The algebra $\slth$ is defined by 
\be 
\slth=\left( \mathfrak{sl}_{2} \otimes \mathbb{C}[t, t^{-1}] \right) 
\oplus \mathbb{C} K \oplus \mathbb{C} D. 
\en 
Here $K$ is the central element. 
We denote the subalgebra 
$\left( \mathfrak{sl}_{2} \otimes \mathbb{C}[t, t^{-1}] \right) 
\oplus \mathbb{C} K $ 
by $\slth'$. 
Set $X(m):=X \otimes t^{m}$ for 
$X=E, F$ or $H$. 
The commutation relations are given by 
\be 
&& 
[ X(m), Y(n)]=[X, Y](m+n)+m\delta_{m+n, 0}(X|Y)K, \\ 
&& 
[ D, X(m) ]=mX(m),  
\en 
where $(X|Y)={\rm tr}(XY)$. 
 
Let $V^{(j)}$ be the $(j+1)$-dimensional irreducible module 
of $\mathfrak{sl}_{2}$. 
We denote its affinization by $V_{z}^{(j)}:=V^{(j)} \otimes \mathbb{C}[z, z^{-1}]$. 
%Then the action of $\slth'$ on $V_{z}^{(j)}$ can be defined by 
%$X(m)v=z^{m}(Xv)$ and $Kv=0$ for $v \in V^{(j)}$. 

In the following, we fix an integer $k \ge 0$ and consider 
highest weight representations of level $k$.  
Let $\lambda$ be an integer such that $0 \le \lambda \le k$. 
We denote the integrable irreducible highest weight module of $\slth$ with 
the highest weight $(k-\lambda)\Lambda_{0}+\lambda \Lambda_{1}$ by $V(\lambda)$. 
Let $\ket{\lambda}$ be the highest weight vector.  

The Virasoro algebra acts on $V(\lambda)$ by the Sugawara operator: 
\be 
L_{n}:=\frac{1}{2(k+2)}\sum_{m \in \mathbb{Z}} 
: E(n-m)F(m)+F(n-m)E(m)+\frac{1}{2}H(n-m)H(m):. 
\en 
Here $: \, \cdot \, :$ is the normal ordering defined by 
\be 
:X(m)Y(n)\! : \, = \, \left\{ 
\begin{array}{ll}  
X(m)Y(n) & (m<n), \\ 
\frac{1}{2}(X(m)Y(m)+Y(m)X(m)) & (m=n), \\ 
Y(n)X(m) & (m>n). 
\end{array} \right. 
\en 
The central charge is given by $c=\frac{3k}{k+2}$. 
The vector $\ket{\lambda}$ satisfies 
\be 
L_{0}\ket{\lambda}=\Delta_{\lambda}\ket{\lambda}, \quad
\Delta_{\lambda}:=\frac{\lambda(\lambda+2)}{4(k+2)}. 
\en 
The module $V(\lambda)$ is bi-graded by $L_{0}$ and $H(0)$. 
We set 
\be 
V(\lambda)_{d, s}:=\{ v \in V(\lambda) \, | \,  L_{0}=dv, \,\, H(0)v=sv  \}. 
\en 
Then $V(\lambda)_{\Delta_{\lambda}, \lambda}=\mathbb{C}\ket{\lambda}$ and 
$V(\lambda)=
\oplus_{d \in\Delta_{\lambda}+\mathbb{Z}_{\ge 0}, \, s \in \lambda+2\mathbb{Z}}V(\lambda)_{d, s}$. 

Consider the dual module 
$V(\lambda)^{*}=\oplus_{d, s}(V(\lambda)_{d, s})^{*}$. 
%For $\bra{u} \in V(\lambda)^{*}$ and $\ket{v} \in V(\lambda)$ 
%we write the dual coupling $V(\lambda)^* \times V(\lambda) \to \mathbb{C}$ as 
%$(\bra{u}, \ket{v}) \mapsto \langle u|v \rangle$. 
%Then $V(\lambda)^{*}$ is the right $\slth$-module with the action 
%$(\bra{u}x)\ket{v}:=\bra{u}(x\ket{v})$. 
It is the right $\slth$-module. 
%and the lowest weight module with the lowest weight $\lambda$. 
Let $\bra{\lambda} \in (V(\lambda)_{\Delta_{\lambda}, \lambda})^{*}$ be the vector 
satisfying $\langle \lambda|\lambda \rangle=1$.

\subsection{Vertex operators} 
Set 
\be 
\mathbb{C}((z)):=\{ \sum_{n \in \mathbb{Z}} a_{n}z^{n} \, | \, a_{n} \in \mathbb{C}, \,\, 
a_{n}=0 \,\, \hbox{for} \,\, n\ll 0 \}. 
\en 
The vertex operator $\phi(z)$ is a $\mathbb{C}[z, z^{-1}]$-linear map 
\be 
\phi(z)=\sum_{n \in \Delta_{\lambda}-\Delta_{\mu}+\mathbb{Z}}
\phi_{n}z^{-n-\Delta_{j}}: 
V(\lambda) \otimes V_{z}^{(j)} \longrightarrow 
V(\mu)\otimes z^{\Delta_{\mu}-\Delta_{\lambda}-\Delta_{j}} \mathbb{C}((z))
\label{eq:def-VO} 
\en 
which commutes with the action of $\slth'$ and satisfies 
\bea 
[L_{n}, \phi(z)]=z^{n}\left(
z\partial+(n+1)\Delta_{j} \right)\phi(z). 
\label{eq:virasoro-commrel}
\ena 
For $u \in V^{(j)}_{z}$ we define the map 
\be 
\phi(u; z):V(\lambda) \longrightarrow 
V(\mu) \otimes z^{\Delta_{\mu}-\Delta_{\lambda}-\Delta_{j}}\mathbb{C}((z)), \quad 
\phi(u; z)v:=\phi(z)(v \otimes u). 
\en 
Then it satisfies 
\bea
[x, \phi(u; z)]=\phi(xu; z) \quad \hbox{for} \quad x \in \slth'. 
\label{eq:phi-commrel} 
\ena 

Consider the function 
\be 
\bra{\mu} \phi(z_{1}) \phi(z_{2}) \ket{\lambda}:
V_{z_{2}}^{(j_{2})} \otimes V_{z_{1}}^{(j_{1})}  
\longrightarrow \mathbb{C}. 
\en 
Note that the relation \eqref{eq:virasoro-commrel} with $n=0$ 
implies 
\be 
\phi(tz)=t^{-\Delta_{j}}\cdot 
t^{L_{0}} \phi(z) t^{-L_{0}}. 
\en 
Hence we have 
\be 
\bra{\mu} \phi(z_{1}) \phi(z_{2}) \ket{\lambda}= 
z_{1}^{\Delta_{\mu}-\Delta_{\lambda}-\Delta_{j_1}-\Delta_{j_{2}}} 
\bra{\mu} \phi(1) \phi(z_{2}/z_{1}) \ket{\lambda}. 
\en 
The function 
\be 
G(z):=\bra{\mu} \phi(1) \phi(z) \ket{\lambda} 
\en 
satisfies the following differential equation: 
\bea 
\frac{1}{\kappa}\frac{dG}{dz}=
\left( \frac{\Omega_{0}}{z}+\frac{\Omega_{1}}{z-1} \right) G, 
\label{eq:KZ} 
\ena 
where 
\be 
&& 
\kappa=\frac{1}{k+2}, \quad 
\Omega_{0}=-F\otimes E-1\otimes FE+\frac{\lambda}{2}(1 \otimes H), \\ 
&& 
\Omega_{1}=E \otimes F+F \otimes E+\frac{1}{2}H \otimes H. 
\en 
%The equation \eqref{eq:KZ} can be derived from the 
%Knizhnik-Zamolodchikov equation. 

Denote by $\phi^{\pm}(z)$ the vertex operator of the following type: 
\be 
\phi^{\pm}(z): V(\lambda) \otimes V_{z}^{(1)} \longrightarrow 
V(\lambda \pm 1) \otimes z^{\Delta_{\lambda \pm 1}-\Delta_{\lambda}-\Delta_{1}} \mathbb{C}((z)). 
\en 
We abbreviate $V^{(1)}$ and $V_{z}^{(1)}$ to $V$ and $V_{z}$, respectively.  
Let $\{ v_{+}, v_{-} \}$ be a basis of $V$ satisfying $Hv_{\pm}=\pm v_{\pm}$. 
We set 
\be 
\phi_{\epsilon}^{\sigma}(z):=\phi^{\sigma}(v_{\epsilon}; z)=
\sum_{n \in \Delta_{\lambda}-\Delta_{\lambda+\sigma}+\mathbb{Z}}
\phi_{\epsilon, n}^{\sigma}z^{-n-\Delta_{1}} 
\en 
for $\sigma, \epsilon=\pm$. 
Each Fourier component of $\phi_{\epsilon}^{\sigma}(z)$ gives  
a map 
\be 
\phi_{\epsilon, n}^{\sigma}: V(\lambda)_{d, s} \longrightarrow 
V(\lambda+\sigma)_{d-n, s+\epsilon}. 
\en  
We choose the normalization 
\bea 
\phi_{\pm, \Delta_{\lambda}-\Delta_{\lambda \pm 1}}^{\pm} \ket{\lambda}=\ket{\lambda \pm 1}. 
\label{eq:normalization++} 
\ena 
Since the vector $\ket{\lambda}$ satisfies 
\be 
F(0)^{\lambda+1}\ket{\lambda}=0, \quad 
E(-1)^{k-\lambda+1}\ket{\lambda}=0, 
\en 
we have 
\bea 
&& 
\phi_{-}^{+}(z)\ket{\lambda}=
z^{\Delta_{\lambda+1}-\Delta_{\lambda}-\Delta_{1}}\left( 
\frac{1}{\lambda+1}F(0)\ket{\lambda+1}+O(z) \right), 
\label{eq:normalization+-}
\\ 
&& 
\phi_{+}^{-}(z)\ket{\lambda}=
z^{\Delta_{\lambda-1}-\Delta_{\lambda}-\Delta_{1}+1}\left( 
\frac{1}{k-\lambda+1}E(-1)\ket{\lambda-1}+O(z) \right). \nn 
\ena 

%{}From the consideration above, 
We can obtain the following formulae by solving \eqref{eq:KZ}. 
Each solution is uniquely determined by \eqref{eq:normalization++} 
and \eqref{eq:normalization+-}. 
\begin{prop}\label{prop:correlation}
\be 
%%%%%%%
\bra{\lambda\pm 2} 
\phi_{\pm}^{\pm}(z_{1}) \phi_{\pm}^{\pm}(z_{2}) \ket{\lambda}&=& 
(z_{1}z_{2})^{\Delta_{\lambda \pm 1}-\Delta_{\lambda}-\Delta_{1}}
(z_{1}-z_{2})^{\frac{\kappa}{2}}, \\  
%%%%%%%
\bra{\lambda} \phi_{+}^{+}(z_{1})\phi_{-}^{-}(z_{2}) \ket{\lambda}&=&
(z_{1}z_{2})^{-\Delta_{1}}
z^{-\frac{1}{2}\left( \lambda+\frac{1}{2} \right) \kappa} 
(1-z)^{-\frac{3\kappa}{2}} \\ 
&& {}\times 
F(-\kappa, 1-(\lambda+2)\kappa, 1-(\lambda+1)\kappa; z), \\ 
%%%%%% 
\bra{\lambda} \phi_{+}^{-}(z_{1})\phi_{-}^{+}(z_{2}) \ket{\lambda}&=&
(z_{1}z_{2})^{-\Delta_{1}}
\frac{-1}{\lambda+1} 
z^{\frac{1}{2}\left( \lambda+\frac{3}{2} \right) \kappa} 
(1-z)^{-\frac{3\kappa}{2}} \\ 
&& {}\times 
F(1-\kappa, \lambda\kappa, 1+(\lambda+1)\kappa; z), \\
%%%%%% 
\bra{\lambda} \phi_{-}^{+}(z_{1})\phi_{+}^{-}(z_{2}) \ket{\lambda}&=&
(z_{1}z_{2})^{-\Delta_{1}}
\frac{-1}{k-\lambda+1} 
z^{1-\frac{1}{2}\left( \lambda+\frac{1}{2} \right) \kappa} 
(1-z)^{-\frac{3\kappa}{2}} \\ 
&& {}\times 
F(1-\kappa, 1-(\lambda+2)\kappa, 2-(\lambda+1)\kappa; z), \\
%%%%%%
\bra{\lambda} \phi_{-}^{-}(z_{1})\phi_{+}^{+}(z_{2}) \ket{\lambda}&=&
(z_{1}z_{2})^{-\Delta_{1}}
z^{\frac{1}{2}\left( \lambda+\frac{3}{2} \right) \kappa} 
(1-z)^{-\frac{3\kappa}{2}} 
F(-\kappa, \lambda\kappa, (\lambda+1)\kappa; z). 
\en 
Here $z=z_{2}/z_{1}$ and $F(a, b, c; z)$ is the hypergeometric function.
\end{prop}

\subsection{Exchange relations} 
In this subsection we prove the exchange relations 
for the vertex operators 
$\phi^{\sigma}_{\epsilon}(z) \, (\sigma, \epsilon=\pm)$. 
To this aim we prepare the following lemma 
which plays a similar role to Lemma \ref{lem:bilinear}. 
\begin{lem}
Let $F(a, b, c; z)$ be the hypergeometric function. 
Set 
\be 
&& 
A(z)=\left( 
\begin{array}{cc} 
\scriptstyle(\beta-\alpha+1)F(\alpha-1, \alpha+\beta-1, 1+\beta; z) & 
\scriptstyle (1-\alpha-\beta)F(\alpha-1, \alpha-\beta-1, 1-\beta; z) \\ 
\scriptstyle F(\alpha, \alpha+\beta, 1+\beta; z) & 
\scriptstyle F(\alpha, \alpha-\beta, 1-\beta; z) 
\end{array} \right)
\en 
and 
\be 
B(z)=\left( 
\begin{array}{cc} 
\scriptstyle\frac{\beta-\alpha}{\beta}F(-\alpha, 1-\alpha-\beta, 1-\beta; z) & 
\scriptstyle \frac{\alpha(\beta-\alpha)}{\beta(1-\beta)}z
F(1-\alpha, 1-\alpha-\beta, 2-\beta; z) \\ 
\scriptstyle\frac{\alpha}{\beta}F(1-\alpha, \beta-\alpha, 1+\beta; z) & 
\scriptstyle F(-\alpha, \beta-\alpha, \beta; z) 
\end{array} \right). 
\en 
Then $A(z)B(z)$ is a polynomial in $z$ of degree one. 
More explicitly we have 
\bea 
A(z)B(z)=\left( 
\begin{array}{cc} 
(1+\beta-3\alpha)-(\alpha+\beta-1)z & 
(1+\beta-3\alpha)z-(\alpha+\beta-1) \\ 
1 & 1 
\end{array} \right). 
\label{eq:hyp-bilinear} 
\ena 
\end{lem} 
\begin{proof} 
By using the formula  
\be 
%%%%%%%%
F(a, b, c; z)&=&
\frac{\Gamma(c)\Gamma(b-a)}{\Gamma(b)\Gamma(c-a)}
(-z)^{-a}F(a, a-c+1, a-b+1; z^{-1}) \\ 
%%%%%%%%
&& {}+ 
\frac{\Gamma(c)\Gamma(a-b)}{\Gamma(a)\Gamma(c-b)}
(-z)^{-b}F(b, b-c+1, b-a+1; z^{-1}) \quad 
(|z|>1, \, z \not\in \mathbb{R}_{>0}), 
\en 
we have 
\be 
A(z)=\left(
\begin{array}{cc} 
z & 0 \\ 0 & 1 
\end{array} \right) 
A(z^{-1})M(\alpha, \beta; z), \quad 
B(z)=
M(\alpha, \beta; z)^{-1}B(z^{-1})
\left(
\begin{array}{cc} 
0 & 1 \\ 1 & 0 
\end{array} \right), 
\en 
where 
\be 
M(\alpha, \beta; z)=\left( 
\begin{array}{cc} 
\displaystyle 
-(-z)^{\alpha+\beta}\frac{\sin{\pi \alpha}}{\sin{\pi \beta}} & 
\displaystyle 
(-z)^{\alpha}
\frac{\Gamma(1-\beta)\Gamma(-\beta)}{\Gamma(1-\alpha-\beta)\Gamma(\alpha-\beta)} \\ 
\displaystyle 
(-z)^{\alpha}
\frac{\Gamma(\beta)\Gamma(1+\beta)}{\Gamma(\beta-\alpha+1)\Gamma(\beta+\alpha)} & 
\displaystyle 
(-z)^{\alpha-\beta}\frac{\sin{\pi \alpha}}{\sin{\pi \beta}} 
\end{array} \right). 
\en 
Hence we obtain 
\be 
A(z)B(z)=
\left(
\begin{array}{cc} 
z & 0 \\ 0 & 1 
\end{array} \right)A(z^{-1})B(z^{-1}) 
\left(
\begin{array}{cc} 
0 & 1 \\ 1 & 0 
\end{array} \right). 
\en 
This equality implies that 
$A(z)B(z)$ is a polynomial of degree one. 
\end{proof} 

\begin{prop}\label{prop:KZquadrel}
We have the following identities of formal power series of 
operators acting on the irreducible highest weight module 
$V(\lambda)$ of level $k$: 
\bea
&& 
h_{\lambda}^{\sigma}(z_{1}, z_{2})
\phi_{\epsilon_{1}}^{\sigma}(z_{1})\phi_{\epsilon_{2}}^{\sigma}(z_{2})-
h_{\lambda}^{\sigma}(z_{2}, z_{1})
\phi_{\epsilon_{2}}^{\sigma}(z_{2})\phi_{\epsilon_{1}}^{\sigma}(z_{1})=0, \quad (\sigma=\pm) 
\label{eq:++rel} \\ 
&& 
\sum_{\sigma=\pm}f_{\lambda}^{\sigma}(z_{1}, z_{2})
\phi_{\epsilon_{1}}^{\sigma}(z_{1})\phi_{\epsilon_{2}}^{-\sigma}(z_{2})-
\sum_{\sigma=\pm}f_{\lambda}^{\sigma}(z_{2}, z_{1})
\phi_{\epsilon_{2}}^{\sigma}(z_{2})\phi_{\epsilon_{1}}^{-\sigma}(z_{1}) 
\label{eq:+-rel1} \\ 
&& {}=
\epsilon_{1}\delta_{\epsilon_{1}+\epsilon_{2}, 0}\delta(z_{2}/z_{1}), \nn \\ 
&& 
\sum_{\sigma=\pm}g_{\lambda}^{\sigma}(z_{1}, z_{2})
\phi_{\epsilon_{1}}^{\sigma}(z_{1})\phi_{\epsilon_{2}}^{-\sigma}(z_{2})-
\sum_{\sigma=\pm}g_{\lambda}^{\sigma}(z_{2}, z_{1})
\phi_{\epsilon_{2}}^{\sigma}(z_{2})\phi_{\epsilon_{1}}^{-\sigma}(z_{1}) 
\label{eq:+-rel2} \\ 
&& {}=
\epsilon_{1}\delta_{\epsilon_{1}+\epsilon_{2}, 0}z_{1}^{-1}\delta(z_{2}/z_{1}). \nn
\ena  
Here $\epsilon_{1}, \epsilon_{2}=\pm, \delta(z_{2}/z_{1})=\sum_{n \in \mathbb{Z}}z_{1}^{-n}z_{2}^{n}$ and 
\be 
%%%%%%% 
h_{\lambda}^{\pm}(z_{1}, z_{2})&=&
(z_{1}z_{2})^{\frac{\kappa}{2}\mp \frac{1}{2}(\lambda+1)\kappa}
z_{1}^{-\frac{\kappa}{2}} \cdot (1-z)^{-\frac{\kappa}{2}}, \\ 
%%%%%%%
f_{\lambda}^{\pm}(z_{1}, z_{2})&=& 
C_{\lambda}^{\pm}(z_{1}z_{2})^{\Delta_{1}} \\ 
&& {}\times 
z^{-\frac{\kappa}{4}\pm\frac{1}{2}(\lambda+1)\kappa}(1-z)^{-1+\frac{3\kappa}{2}} 
F(-1+\kappa, -1+\kappa\pm(\lambda+1)\kappa, 1\pm(\lambda+1)\kappa; z), \\ 
%%%%%%%
g_{\lambda}^{\pm}(z_{1}, z_{2})&=& 
\tilde{C}_{\lambda}^{\pm}z_{1}^{-1}(z_{1}z_{2})^{\Delta_{1}} \\ 
&& {}\times 
z^{-\frac{\kappa}{4}\pm\frac{1}{2}(\lambda+1)\kappa}(1-z)^{-1+\frac{3\kappa}{2}}
F(\kappa, \kappa\pm(\lambda+1)\kappa, 1\pm(\lambda+1)\kappa; z), 
\en
where $z=z_{2}/z_{1}$ and $F(a, b, c; z)$ is the hypergeometric function. 
The constants $C_{\lambda}^{\pm}, \tilde{C}_{\lambda}^{\pm}$ are given by 
\be 
C_{\lambda}^{+}=\frac{\lambda(k+\lambda+2)}{2k(\lambda+1)}, \quad 
C_{\lambda}^{-}=\frac{\lambda-k}{2k}, \quad 
\tilde{C}_{\lambda}^{+}=\frac{\lambda}{\lambda+1}, \quad 
\tilde{C}_{\lambda}^{-}=-1. 
\en 
\end{prop}

\begin{rem} 
In the case of $\lambda=0$, the summation in the left hand side of 
\eqref{eq:+-rel1} and \eqref{eq:+-rel2} becomes one term with $\sigma=-$. 
Similarly, if $\lambda=k$, it becomes one term with $\sigma=+$. 
\end{rem} 
 
\begin{proof}[\normalfont\bfseries Proof of Proposition \ref{prop:KZquadrel}]
Here we prove the second relation \eqref{eq:+-rel1}. 
The proofs for the other ones are similar. 

Let $\langle \cdot , \cdot \rangle$ be the bilinear pairing 
$V \times V \to \mathbb{C}$ defined by 
$\langle v_{\epsilon_{1}}, v_{\epsilon_{2}} \rangle:=
 \epsilon_{1}\delta_{\epsilon_{1}+\epsilon_{2}, 0}$. 
This pairing satisfies 
\bea 
\langle Xu, v \rangle+\langle u, Xv \rangle=0, \quad 
u, v \in V 
\label{eq:pairing-invariance} 
\ena 
for $X=E, F, H$.  

For $u_{i} \in V \, (i=1, 2)$ we set 
\be 
%%%%%%%
A_{u_{1}, u_{2}}(z_{1}, z_{2})&:=&
\sum_{\sigma=\pm}f_{\lambda}^{\sigma}(z_{1}, z_{2})
\phi^{\sigma}(u_{1}; z_{1})\phi^{-\sigma}(u_{2}; z_{2}) \\ 
%%%%%%% 
&&{}-
\sum_{\sigma=\pm}f_{\lambda}^{\sigma}(z_{2}, z_{1})
\phi^{\sigma}(u_{2}; z_{2})\phi^{-\sigma}(u_{1}; z_{1})-
\langle u_{1}, u_{2} \rangle \delta(z_{2}/z_{1}). 
\en 
Then \eqref{eq:+-rel1} is equivalent to the equality 
$A_{u_{1}, u_{2}}(z_{1}, z_{2})=0$. 
Let us prove it. 
We show the following identities: 
\bea 
&& 
\bra{\lambda} A_{u_{1}, u_{2}}(z_{1}, z_{2}) \ket{\lambda}=0, 
\label{eq:A1} \\ 
&& 
[X(m), A_{u_{1}, u_{2}}(z_{1}, z_{2})]=
z_{1}^{m}A_{Xu_{1}, u_{2}}(z_{1}, z_{2})+z_{2}^{m}A_{u_{1}, Xu_{2}}(z_{1}, z_{2}), \quad 
\hbox{for} \,\, X=E, F, H. 
\label{eq:A2}
\ena 
These identities imply that $A_{u_{1}, u_{2}}(z_{1}, z_{2})=0$ 
by the same argument as the proof of Proposition \ref{QUADREL}. 

To prove \eqref{eq:A1} it is enough to consider the case 
that $u_{1}=v_{\epsilon}$ and $u_{2}=v_{-\epsilon}$ for $\epsilon=\pm$. 
Then \eqref{eq:A1} follows from Proposition \ref{prop:correlation} 
and the equalities of  
(1, 1) and (1, 2) elements in \eqref{eq:hyp-bilinear} 
with $\alpha=\kappa$ and $\beta=(\lambda+1)\kappa$. 

Let us prove \eqref{eq:A2}. 
{}From \eqref{eq:phi-commrel} we have 
\be 
%%%%%%%
[X(m), A_{u_{1}, u_{2}}(z_{1}, z_{2})]&=&
z_{1}^{m}\left( 
A_{Xu_{1}, u_{2}}(z_{1}, z_{2})+\langle Xu_{1}, u_{2} \rangle \delta(z_{2}/z_{1}) \right) \\ 
&& {}+z_{2}^{m}\left( 
A_{u_{1}, Xu_{2}}(z_{1}, z_{2})+\langle u_{1}, Xu_{2} \rangle \delta(z_{2}/z_{1}) \right). 
\en 
Note that $z_{1}^{m}\delta(z_{2}/z_{1})=z_{2}^{m}\delta(z_{2}/z_{1})$. 
Hence \eqref{eq:A2} follows form \eqref{eq:pairing-invariance}. 
This completes the proof of \eqref{eq:+-rel1}. 

In the proof of \eqref{eq:+-rel2} we should show that 
\be 
&& 
\bra{\lambda}\left( 
\sum_{\sigma=\pm}g_{\lambda}^{\sigma}(z_{1}, z_{2})
\phi_{\epsilon}^{\sigma}(z_{1})\phi_{-\epsilon}^{-\sigma}(z_{2})-
\sum_{\sigma=\pm}g_{\lambda}^{\sigma}(z_{2}, z_{1})
\phi_{-\epsilon}^{\sigma}(z_{2})\phi_{\epsilon}^{-\sigma}(z_{1})
\right) \ket{\lambda} \\ 
&& {}=
\epsilon_{1}\delta_{\epsilon_{1}+\epsilon_{2}, 0}z_{1}^{-1}\delta(z_{2}/z_{1}). 
\en 
To prove this, use the equalities of (2, 1) and (2, 2) elements 
in \eqref{eq:hyp-bilinear}. 
\end{proof}

\subsection{Monomial basis} 
{}From the exchange relations \eqref{eq:++rel}, \eqref{eq:+-rel1} and 
\eqref{eq:+-rel2}, we can construct a spanning set of $V(\mu)$: 

\begin{prop} 
The vectors 
\bea 
\phi_{\epsilon_{1}, n_{1}}^{\sigma_{1}} \cdots 
\phi_{\epsilon_{L}, n_{L}}^{\sigma_{L}} \ket{0}, \quad 
L \ge 0, \,\, \sigma_{i} \in \{ +, - \}, \,\, \epsilon_{i} \in \{ +, -\}
\label{eq:KZmonomial-basis} 
\ena 
satisfying the following condition \eqref{eq:path-condition} and 
\eqref{eq:weight-condition} span the irreducible highest weight module $V(\mu)$ of level $k$. 
Set $\lambda_{L}=0$ and $\lambda_{i-1}=\lambda_{i}+\sigma_{i} \,\, (i=1, \ldots , L)$. 
Then 
\bea 
&& 
0 \le \lambda_{i} \le k, \quad 
\lambda_{0}=\mu, \quad n_{i} \in \Delta_{\lambda_{i}}-\Delta_{\lambda_{i-1}}+\mathbb{Z}, 
\label{eq:path-condition} \\ 
&& 
-n_{L} \ge \Delta_{1}, \quad 
n_{i+1}-n_{i} \ge w(\lambda_{i-1}, \lambda_{i}, \lambda_{i+1})+h(\epsilon_{i}, \epsilon_{i+1}) 
\,\, (i=1, \ldots , L-1), 
\label{eq:weight-condition} 
\ena 
where the functions $w(\lambda_{i-1}, \lambda_{i}, \lambda_{i+1})$ and 
$h(\epsilon_{i}, \epsilon_{i+1})$ are defined by 
\be 
&& 
w(\lambda\pm1, \lambda, \lambda\mp1)=\frac{\kappa}{2}, \quad 
w(\lambda, \lambda-1, \lambda)=( \lambda+\frac{1}{2} ) \kappa, \quad 
w(\lambda, \lambda+1, \lambda)=1-( \lambda+\frac{3}{2})\kappa, \\ 
&& 
h(+, -)=1, \quad h(\epsilon, \epsilon')=0 \,\, \hbox{otherwise}. 
\en 
\end{prop} 
\begin{proof} 
First we prove that the vectors 
$\phi_{\epsilon_{1}, n_{1}}^{\sigma_{1}} \cdots 
\phi_{\epsilon_{L}, n_{L}}^{\sigma_{L}} \ket{0}$ 
satisfying only the condition \eqref{eq:path-condition} span $V(\mu)$. 
Let $W_{\mu}$ be the subspace of $V(\mu)$ spanned by the vectors. 
{}From \eqref{eq:normalization++} we have $\ket{\mu} \in W_{\mu}$. 
Moreover, from the commutation relation \eqref{eq:phi-commrel} and 
\be 
E(0)\ket{0}=F(0)\ket{0}=F(1)\ket{0}=0, \quad 
E(-1)\ket{0}=k\phi_{+, \Delta_{1}-1}^{-}\phi_{+, -\Delta_{1}}^{+}\ket{0}, 
\en 
which follows from \eqref{eq:normalization+-}, 
we see that $xW_{\mu} \subset W_{\mu}$ for any $x \in \slth$. 
Hence $W_{\mu}=V(\mu)$.  

The rest of the proof is similar to that in Section \ref{sec:exchange}. 
We consider quadratic monomials 
$\phi_{\epsilon_{1}, n_{1}}^{\sigma_{1}}\phi_{\epsilon_{2}, n_{2}}^{\sigma_{2}}$
acting on $V(\lambda)$, and show that each monomial is reduced to 
a linear combination of ones satisfying 
\be 
n_{2}-n_{1} 
\ge w(\lambda+\sigma_{1}+\sigma_{2}, \lambda+\sigma_{2}, \lambda)+h(\epsilon_{1}, \epsilon_{2}) 
\en 
by using the exchange relations in Proposition \ref{prop:KZquadrel}. 
Here we prove it in the case of $\sigma_{1}=-\sigma_{2}$ and $0<\lambda<k$. 
The proofs for the other cases are similar. 

Fix $n_{1}+n_{2}$ and set 
\be 
\Phi_{\epsilon_{1}, \epsilon_{2}; n_{2}-n_{1}}^{\sigma}=
\phi_{\epsilon_{1}, n_{1}}^{\sigma}\phi_{\epsilon_{2}, n_{2}}^{-\sigma}. 
\en 
%We say $\Phi_{\epsilon_{1}, \epsilon_{2}; n}^{\sigma}$ is larger than 
%$\Phi_{\epsilon_{1}', \epsilon_{2}'; n'}^{\sigma'}$ if $n<n'$. 
We say the monomial $\Phi_{\epsilon_{1}, \epsilon_{2}; n}^{\sigma}$ is admissible if 
\be 
n \ge w(\lambda, \lambda-\sigma, \lambda)+h(\epsilon_{1}, \epsilon_{2}). 
\en 

Write down the exchange relations \eqref{eq:+-rel1} and \eqref{eq:+-rel2} 
in terms of Fourier components. Then we find 
\bea
&& 
C_{\lambda}^{+}\Phi_{\epsilon_{1},\epsilon_{2}; n}^{+}+
C_{\lambda}^{-}\Phi_{\epsilon_{1},\epsilon_{2}; n-2(\lambda+1)\kappa}^{-}+\cdots  
\label{eq:phin-quadrel1} \\ 
&& {}-\left(
C_{\lambda}^{+}\Phi_{\epsilon_{2}, \epsilon_{1}; -n+(2\lambda+1)\kappa}^{+}+
C_{\lambda}^{-}\Phi_{\epsilon_{2}, \epsilon_{1}; -n-\kappa}^{-}+ \cdots 
\right)\equiv 0 \nn
%&& {}=
%\epsilon_{1}\delta_{\epsilon_{1}+\epsilon_{2}}\delta_{m+n, 0}, \nn
\ena 
and 
\bea 
&& 
\tilde{C}_{\lambda}^{+}\Phi_{\epsilon_{1}, \epsilon_{2}; n}^{+}+
\tilde{C}_{\lambda}^{-}\Phi_{\epsilon_{1}, \epsilon_{2}; n-2(\lambda+1)\kappa}^{-}+\cdots 
\label{eq:phin-quadrel2} \\ 
&& {}-\left( 
\tilde{C}_{\lambda}^{+}\Phi_{\epsilon_{2}, \epsilon_{1}; -n+(2\lambda+1)\kappa+2}^{+}+
\tilde{C}_{\lambda}^{-}\Phi_{\epsilon_{2}, \epsilon_{1}; -n-\kappa+2}^{-}+\cdots 
\right)\equiv 0. \nn
%&& {}=
%\epsilon_{1}\delta_{\epsilon_{1}+\epsilon_{2}}\delta_{m+n, 0}. \nn
\ena 
Here $\equiv$ means that we ignore the constant terms. 
{}From these relations and 
\bea 
\det{ \left( 
\begin{array}{cc} 
C_{\lambda}^{+} & C_{\lambda}^{-} \\ \tilde{C}_{\lambda}^{+} & \tilde{C}_{\lambda}^{-} 
\end{array} 
\right) }=-\frac{\lambda}{k} \not=0, 
\label{eq:det-nonzero} 
\ena 
the monomials 
\be 
&& 
\Phi_{\epsilon_{1}, \epsilon_{2}; n}^{+}, \quad n<(\lambda+\frac{1}{2})\kappa=
w(\lambda, \lambda-1, \lambda) 
\quad \hbox{and} \\ 
&& 
\Phi_{\epsilon_{1}, \epsilon_{2}; n}^{-}, \quad n<-(\lambda+\frac{3}{2})\kappa=
w(\lambda, \lambda+1, \lambda)-1 
\en 
can be written as a linear combination of the rest. 
Hence it suffices to prove that the monomials 
\be 
\Phi_{\epsilon_{1}, \epsilon_{2}; -(\lambda+3/2)\kappa}^{-}, \quad 
\Phi_{+-; (\lambda+1/2)\kappa}^{+} \quad 
\hbox{and} \quad  
\Phi_{+-; -(\lambda+3/2)\kappa+1}^{-}
\en 
can be reduced to admissible ones. 

Set $\epsilon_{1}=\epsilon_{2}=\epsilon$ and $n=(\lambda+1/2)\kappa$ 
in \eqref{eq:phin-quadrel2}. 
Then we have 
\be 
\tilde{C}_{\lambda}^{+}\Phi_{\epsilon, \epsilon; (\lambda+1/2)\kappa}^{+}+
\tilde{C}_{\lambda}^{-}\Phi_{\epsilon, \epsilon; -(\lambda+3/2)\kappa}^{-}+\cdots 
\equiv 0. 
\en 
Hence $\Phi_{\pm, \pm; -(\lambda+3/2)\kappa}^{-}$ can be reduced. 

Next consider \eqref{eq:phin-quadrel1} and \eqref{eq:phin-quadrel2} with 
$\epsilon_{1}=-\epsilon_{2}=\epsilon$ and $n=(\lambda+1/2)\kappa$: 
\bea 
&&
C_{\lambda}^{+}\left( 
\Phi_{\epsilon, -\epsilon; (\lambda+1/2)\kappa}^{+}-
\Phi_{-\epsilon, \epsilon; (\lambda+1/2)\kappa}^{+} \right) 
\label{eq:border1} \\
&& {}+
C_{\lambda}^{-}\left( 
\Phi_{\epsilon, -\epsilon; -(\lambda+3/2)\kappa}^{-}- 
\Phi_{-\epsilon, \epsilon; -(\lambda+3/2)\kappa}^{-} \right)+\cdots  
\equiv 0 \nn
\ena 
and 
\be 
\tilde{C}_{\lambda}^{+}\Phi_{\epsilon, -\epsilon; (\lambda+1/2)\kappa}^{+}+
\tilde{C}_{\lambda}^{-}\Phi_{\epsilon, -\epsilon; -(\lambda+3/2)\kappa}^{-}+\cdots 
\equiv 0. 
\en 
These relations hold for $\epsilon=+$ and $-$. 
Here we note that \eqref{eq:border1} with $\epsilon=+$ and the one with $\epsilon=-$ 
are equivalent. 
Hence we have three relations among four monomials 
$\Phi_{\pm, \mp; (\lambda+1/2)\kappa}^{+}$ and $\Phi_{\pm, \mp; -(\lambda+3/2)\kappa}^{-}$. 
The relations are linearly independent from \eqref{eq:det-nonzero}. 
By using them we reduce the monomials except $\Phi_{-, +; (\lambda+1/2)}^{+}$. 

At last we consider \eqref{eq:phin-quadrel2} with $\epsilon_{1}=-\epsilon_{2}=\epsilon$ 
and $n=(\lambda+1/2)\kappa+1$. Then we find 
\be 
\tilde{C}_{\lambda}^{-}\left( 
\Phi_{\epsilon, -\epsilon; -(\lambda+3/2)\kappa+1}^{-}- 
\Phi_{-\epsilon, \epsilon; -(\lambda+3/2)\kappa+1}^{-}\right) 
+\cdots \equiv 0. 
\en
This relation holds for $\epsilon=+$ and $-$; however 
this gives only one relation by the same reason as before. 
Hence we can reduce one monomial. 
Here we reduce $\Phi_{+-; -(\lambda+3/2)\kappa+1}^{-}$. 
This completes the proof. 
\end{proof}

\begin{thm}\label{thm:KZfinal}
The vectors \eqref{eq:KZmonomial-basis} 
satisfying \eqref{eq:path-condition} and \eqref{eq:weight-condition} 
are linearly independent. 
Hence they give a basis of the irreducible highest weight module $V(\mu)$ 
of level $k$. 
\end{thm}
\begin{proof}  
We consider the character of $V(\mu)$, 
which is a formal power series defined by 
\be 
{\rm ch}_{q, z}V(\mu)=q^{-\Delta_{\mu}}\sum_{d, s}(\dim{V(\mu)_{d, s}})q^{d}z^{s}. 
\en 
%Now let us prove the statement. 
Since the vectors \eqref{eq:KZmonomial-basis} span $V(\mu)$, we have 
\bea 
%%%%%
{\rm ch}_{q, z}V(\mu) &\le& 
%q^{-\Delta_{\mu}}
\delta_{\mu, 0}+q^{-\Delta_{\mu}}\sum_{L \ge 1}\frac{1}{(q)_{L}} 
\sum_{(\lambda_{i})}
q^{L\Delta_{1}+\sum_{i=1}^{L-1}iw(\lambda_{i-1}, \lambda_{i}, \lambda_{i+1})} 
\label{eq:KZcharacter} \\ 
%%%%%
&& \qquad \qquad {}\times 
\sum_{\epsilon_{1}, \ldots , \epsilon_{L}=\pm}
q^{\sum_{i=1}^{L-1}ih(\epsilon_{i}, \epsilon_{i+1})}
z^{\sum_{i=1}^{L}\epsilon_{i}}. \nn
\ena 
Here the second sum is over the sequences 
$(\lambda_{0}, \ldots , \lambda_{L})$ 
of integers such that 
\be 
\lambda_{0}=\mu, \quad \lambda_{i+1}=\lambda_{i}\pm1, \quad 0 \le \lambda_{i} \le k, 
\quad \lambda_{L}=0. 
\en 

Let us prove the equality in \eqref{eq:KZcharacter}. 
Then it implies Theorem \ref{thm:KZfinal}. 
Set 
\be 
\tilde{w}(\lambda_{i-1}, \lambda_{i}, \lambda_{i+1}):=
w(\lambda_{i-1}, \lambda_{i}, \lambda_{i+1})-\Delta_{\lambda_{i-1}}+
2\Delta_{\lambda_{i}}-\Delta_{\lambda_{i+1}}. 
\en 
Then we have 
\be 
\tilde{w}(\lambda_{i-1}, \lambda_{i}, \lambda_{i+1})=\left\{ 
\begin{array}{ll} 
1, & \lambda_{i-1}=\lambda_{i+1}=\lambda_{i}-1, \\ 
0, & \hbox{otherwise}. 
\end{array} 
\right. 
\en 
Hence we find 
\be 
%%%%%%%%%
\sum_{(\lambda_{i})}
q^{L\Delta_{1}+\sum_{i=1}^{L-1}iw(\lambda_{i-1}, \lambda_{i}, \lambda_{i+1})}
&=&
q^{\Delta_{\mu}}\sum_{(\lambda_{i})}
q^{\sum_{i=1}^{L-1}i\tilde{w}(\lambda_{i-1}, \lambda_{i}, \lambda_{i+1})} \\ 
%%%%%%%%%
&=&q^{\Delta_{\mu}}K_{\mu, (1^{L})}^{(k)}(q), 
\en 
where $K_{\mu, (1^{L})}^{(k)}(q)$ is the level-restricted Kostka polynomial. 
{}From this formula and 
\be 
\sum_{\epsilon_{1}, \ldots , \epsilon_{L}=\pm}
q^{\sum_{i=1}^{L-1}ih(\epsilon_{i}, \epsilon_{i+1})}
z^{\sum_{i=1}^{L}\epsilon_{i}}=
\sum_{l=0}^{L}\frac{(q)_{L}}{(q)_{l}(q)_{L-l}}z^{L-2l}, 
\en 
the right hand side of \eqref{eq:KZcharacter} is equal to 
\bea 
%{\rm ch}_{q, z}V(\mu) \le q^{\Delta_{\mu}}
\sum_{L \ge 0}\frac{1}{(q)_{L}}K_{\mu, (1^{L})}^{(k)}(q)
\sum_{l=0}^{L}\frac{(q)_{L}}{(q)_{l}(q)_{L-l}}z^{L-2l}. 
\label{eq:last-formula}
\ena 
Here by definition we set $K_{\mu, (1^{0})}^{(k)}(q)=\delta_{\mu, 0}$. 
Then \eqref{eq:last-formula} is equal to the character of $V(\mu)$ 
as shown in \cite{FJLM} (the formula (2.14) in the limit $N \to \infty$). 
\end{proof}
%%%%%%%%%%%%%%%%%%%%%%%%%%%%%%%%%%%%%%%%%%%%%%%%%%%%%%%%%%%%%%%
%                               %
% Acknowledgments                      %
%                               %
%%%%%%%%%%%%%%%%%%%%%%%%%%%%%%%%%%%%%%%%%%%%%%%%%%%%%%%%%%%%%%%
\bigskip 
\noindent
{\it Acknowledgments.}\quad
BF is partially supported by grants RFBR-02-01-01015,
RFHR-01-01-00906, INTAS-00-00055. JM is partially supported by 
the Grant-in-Aid for Scientific Research (B2) no.12440039, 
and TM is partially supported by 
(A1) no.13304010, Japan Society for the Promotion of Science.
EM is partially supported by the National Science
Foundation (NSF) grant DMS-0140460.

%%%%%%%%%%%%%%%%%%%%%%%%%%%%%%%%%%%%%%%%%%%%%%%%%%%%%%%%%%%%%%%%%%

\end{document}